\documentclass[12pt]{article}
\usepackage{graphicx}

\usepackage{epsfig}
\usepackage{amsmath}
\usepackage{amssymb}
\usepackage{amsthm}
\usepackage{latexsym}
\usepackage{cite}
\usepackage{amsbsy,amsfonts,euscript,eepic,rotating}

\newcommand{\comment}[1]{}
\newtheorem{theorem}{Theorem}
\newtheorem{lemma}{Lemma}[section]
\newtheorem{remark}{Remark}[section]
\newtheorem{corollary}{Corollary}[section]

\newtheorem{definition}{Definition}[section]
\begin{document}

\title{\LARGE
{\bf Two-Dimensional Critical Percolation: \\ The Full Scaling Limit}
}

\author{
{\bf Federico Camia}
\thanks{Research partially supported by a Marie Curie Intra-European Fellowship
under contract MEIF-CT-2003-500740 and by a Veni grant of the Dutch Organization
for Scientific Research (NWO).}\,
\thanks{E-mail: fede@few.vu.nl}\\
{\small \sl Department of Mathematics, Vrije Universiteit Amsterdam}\\
\and
{\bf Charles M.~Newman}
\thanks{Research partially supported by the
U.S. NSF under grant DMS-01-04278.}\,
\thanks{E-mail: newman@courant.nyu.edu}\\
{\small \sl Courant Inst.~of Mathematical Sciences,
New York University}
}

\date{}

\maketitle

\begin{abstract}
We use $SLE_6$ paths to construct a process of continuum nonsimple loops
in the plane and prove that this process coincides with the full continuum
scaling limit of 2D critical site percolation on the triangular lattice --
that is, the scaling limit of the set of all interfaces between different
clusters.
Some properties of the loop process, including conformal invariance,
are also proved.
%
\end{abstract}



\section{Introduction and Main Results} \label{intro}
In the theory of critical phenomena it is usually assumed
that a physical system near a continuous phase transition
is characterized by a single length scale (the ``correlation
length'') in terms of which all other lengths should be measured.
When combined with the experimental observation that the
correlation length diverges at the phase transition, this
simple but strong assumption, known as the scaling hypothesis,
leads to the belief that at criticality the system has
no characteristic length, and is therefore invariant under
scale transformations.
This suggests that all thermodynamic functions at criticality
are homogeneous functions, and predicts the appearance of
power laws.
It also means that it should be possible to rescale a
critical system appropriately and obtain a continuum model
(the ``continuum scaling limit'') which may have more
symmetries and be easier to study than the original discrete
model defined on a lattice.

Indeed, thanks to the work of Polyakov~\cite{polyakov}
and others~\cite{bpz1,bpz2}, it was understood by physicists
since the early seventies that critical statistical mechanical
models should possess continuum scaling limits with a global
conformal invariance that goes beyond simple scale invariance,
as long as the discrete models have ``enough'' rotation invariance.
This property gives important information, enabling the determination
of two- and three-point functions at criticality, when they are
nonvanishing.
Because the conformal group is in general a finite dimensional Lie
group, the resulting constraints are limited in number; however,
the situation becomes particularly interesting in two dimensions,
since there every analytic function $\omega=f(z)$ defines a conformal
transformation, at least at points where $f'(z) \neq 0$.
As a consequence, the conformal group in two dimensions is
infinite-dimensional.

After this observation was made, a large number of critical
problems in two dimensions were analyzed using conformal
methods, which were applied, among others, to Ising and Potts
models, Brownian motion, Self-Avoiding Walk (SAW), percolation,
and Diffusion Limited Aggregation (DLA).
The large body of knowledge and techniques that resulted, starting
with the work of Belavin, Polyakov and Zamolodchikov~\cite{bpz1,bpz2}
in the early eighties, goes under the name of Conformal Field Theory
(CFT).
In two dimensions, one of the main goals of CFT and its most important
application to statistical mechanics is a complete classification
of all universality classes via irreducible representations of the
infinite-dimensional Virasoro algebra.

Partly because of the success of CFT, work in recent
years on critical phenomena seemed to slow down somewhat,
probably due to the feeling that most of the leading
problems had been resolved.
Nonetheless, however powerful and successful it may be,
CFT has some limitations and leaves various open problems.
First of all, the theory deals primarily with correlation
functions of {\it local} (or quasi-local) operators, and is
therefore not always the best tool to investigate other
quantities.
Secondly, given some critical lattice model, there is no
way, within the theory itself, of deciding to which CFT
it corresponds.
A third limitation, of a different nature, is due to the fact
that the methods of CFT, although very powerful, are generally
speaking not completely rigorous from a mathematical point
of view.

In a somewhat surprising twist, the most recent developments
in the area of two-dimensional critical phenomena have
emerged in the mathematics literature and have followed
a new direction, which has provided new tools and a way of
coping with at least some of the limitations of CFT.
The new approach may even provide a reinterpretation of CFT,
and seems to be complementary to the traditional one in the sense
that questions that are difficult to pose and/or answer within
CFT are easy and natural in this new approach and vice versa.

These new developments came on the heels of interesting
results on the scaling limits of discrete models
(see, e.g., the work of Aizenman~\cite{aizenman,aizenman1},
Benjamini-Schramm~\cite{bs}, Aizenman-Burchard~\cite{ab},
Aizenman-Burchard-Newman-Wilson~\cite{abnw},
Aizenman-Duplantier-Aharony~\cite{ada} and
Kenyon~\cite{kenyon1,kenyon2}) but they differ greatly
from those because they are based on a radically new approach
whose main tool is the Stochastic Loewner Evolution ($SLE$),
or Schramm-Loewner Evolution, as it is also known, introduced
by Schramm~\cite{schramm}.
The new approach, which is probabilistic in nature,
focuses directly on non-local structures that characterize
a given system, such as cluster boundaries in Ising, Potts
and percolation models, or loops in the $O(n)$ model.
At criticality, these non-local objects become, in the
continuum limit, random curves whose distributions can be
uniquely identified thanks to their conformal invariance
and a certain ``Markovian" property.
There is a one-parameter family of $SLE$s, indexed by
a positive real number $\kappa$, and they appear to be
the only possible candidates for the scaling limits of
interfaces of two-dimensional critical systems that are
believed to be conformally invariant.

In particular, substantial progress has been made in
recent years, thanks to $SLE$, in understanding the
fractal and conformally invariant nature of (the scaling
limit of) large percolation clusters, which has attracted
much attention and is of interest both for intrinsic
reasons, given the many applications of percolation,
and as a paradigm for the behavior of other systems.
The work of Schramm~\cite{schramm} and Smirnov~\cite{smirnov}
has identified the scaling limit of a certain percolation
interface with $SLE_6$, providing, along with the work of
Lawler-Schramm-Werner~\cite{lsw1,lsw5} and Smirnov-Werner~\cite{sw},
a confirmation of many results in the physics literature,
as well as some new results.

However, $SLE_6$ describes a single interface, which
can be obtained by imposing special boundary conditions,
and is not in itself sufficient to immediately describe
the scaling limit of the unconstrained model (without
boundary conditions) in the whole plane.
In particular, not only the nature and properties, but the
very existence of the scaling limit of the collection of
all interfaces remained an open question.
This is true of all models, such as Ising and Potts models,
that are represented in terms of clusters, and where the set
of all interfaces forms a collection of loops.
As already indicated by Smirnov~\cite{smirnov-long}, such a
collection of loops should have a continuum limit, that
we will call the ``full" scaling limit of the model.
The single interface limit is ideal for analyzing certain
crossing/connectivity probabilities but not so good for others;
in Section~\ref{quickresults} we give a few examples showing the
use of the full scaling limit to represent such probabilities.
In the context of percolation, in~\cite{cn} the authors used
$SLE_6$ to construct a random process of continuous loops in
the plane, which was identified with the full scaling limit
of critical two-dimensional percolation, but without detailed
proofs.
(For a discussion of whether this full scaling limit is a ``black
noise," see~\cite{tsirelson}. For an analysis of random processes
of loops related to $SLE_{\kappa}$ for other values of $\kappa$,
and \emph{conjectured} to correspond to the full scaling limits
of other statistical mechanics models, see~\cite{werner3,ss,shw}.)

In this paper, we complete the analysis of~\cite{cn}, making
rigorous the connection between the construction given there
and the full scaling limit of percolation, and we prove some
properties of the full scaling limit, the Continuum Nonsimple
Loop process, including (one version of) conformal invariance.
The present work, as well as that of Smirnov~\cite{smirnov,smirnov-long},
builds on a collection of papers, including~\cite{lps,aizenman,aizenman1,ksz,ab,ada},
which provided both inspiration and essential technical results.
The proofs are based on the fact that the percolation exploration
path converges in distribution to the trace of chordal $SLE_6$,
as argued by Schramm and
Smirnov~\cite{schramm,smirnov,smirnov-long,smirnov1,smirnov-private-comm},
and in particular on a specific version of this convergence that
we will call statement (S) (see Section~\ref{main-results}).
We note that no detailed proof of any version of convergence
to $SLE_6$ has been available.
Nevertheless, at the request of the editor, we do not include
a detailed proof of statement (S) in the present paper, as
originally planned~\cite{cn1}, due to length considerations.

However, a detailed proof of statement (S), based on Smirnov's
theorem about the convergence of crossing probabilities to Cardy's
formula~\cite{smirnov}, 
is now 
the topic of a separate paper~\cite{cn2}.
We note that statement (S) is restricted to Jordan domains while
no such restriction is indicated in~\cite{smirnov,smirnov-long}.



The rest of the paper is organized as follows.
In Section~\ref{quickresults}, we provide a quick presentation
of Theorems~\ref{main-result1}, \ref{features}, and \ref{main-result2},
which represent most of our main results, with some definitions
postponed until Section~\ref{defs}, including that of $SLE_6$.
Section~\ref{construct} is devoted to the construction mentioned
in Theorem~\ref{main-result2} of the Continuum Nonsimple Loop
process in a finite region $D$ of the plane.
In Section~\ref{lapa}, we introduce the discrete lattice model
and a discrete construction analogous to the continuum one
presented in Section~\ref{construct}.
Most of the main technical results of this paper are stated in
Section~\ref{main-results}, while Section~\ref{proofs} contains
the proofs, using (S), of those results and the
results in Section~\ref{quickresults}.

\subsection{Main Results} \label{quickresults}

\begin{figure}[!ht]
\begin{center}
\includegraphics[width=8cm]{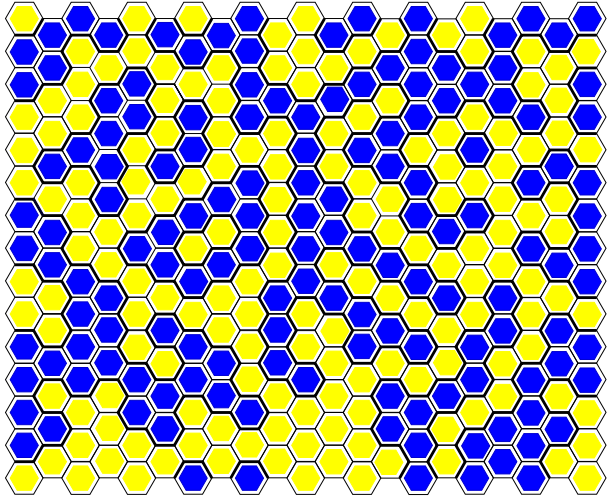}
\caption{
Finite portion of a (site) percolation configuration
on the triangular lattice $\cal T$. Each hexagon of the hexagonal
lattice $\cal H$ represents a site of $\cal T$ and is assigned one
of two colors. In the critical percolation model, colors are assigned
randomly with equal probability.
The cluster boundaries are indicated by heavy lines; some small
loops appear, while other boundaries extend beyond the finite
window.
}
\label{loops}
\end{center}
\end{figure}

At the percolation critical point, with probability one there is
no infinite cluster (in two dimensions), therefore the percolation
cluster boundaries form loops (see Figure~\ref{loops},
where site percolation on the triangular lattice $\cal T$
is depicted exploiting the duality between $\cal T$ and the hexagonal
lattice $\cal H$).
We will refer to the continuum scaling limit (as the mesh size $\delta$
of the rescaled hexagonal lattice $\delta{\cal H}$ goes to zero) of the
collection of all these loops as the Continuum Nonsimple Loop process.
Its existence is the content of Theorem~\ref{main-result1} and some
of its properties are described in Theorem~\ref{features} below.

The Continuum Nonsimple Loop process can be described as a ``gas"
of loops, or more precisely, a probability measure on countable
collections of continuous, nonsimple, fractal loops in the plane.
Later in this paper,
we will provide precise definitions of the objects involved in
the next three theorems as well as detailed proofs.

\begin{theorem} \label{main-result1}
In the continuum scaling limit, the probability distribution of the
collection of all boundary contours of critical site percolation on
the triangular lattice converges to a probability distribution on
collections of continuous, nonsimple loops.
\end{theorem}

\begin{theorem} \label{features}
The Continuum Nonsimple Loop process whose distribution is
specified in Theorem~\ref{main-result1} has the following
properties, which are valid with probability one:
\begin{enumerate}
\item It is a random collection of countably many
noncrossing continuous loops in the plane.
The loops can and do touch themselves and each other
many times, but there are no triple points; i.e. no three
or more loops can come together at the same point, and a
single loop cannot touch the same point more than twice,
nor can a loop touch a point where another loop touches
itself.
\item Any deterministic point $z$ in the plane (i.e.,
chosen independently of the loop process) is surrounded by
an infinite family of nested loops with diameters going to
both zero and infinity; any annulus about that point with
inner radius $r_1 > 0$ and outer radius $r_2 < \infty$
contains only a finite number $N(z, r_1, r_2)$ of those loops.
Consequently, any two distinct deterministic points of the
plane are separated by loops winding around each of them.
\item Any two loops are connected by a finite ``path''
of touching loops.
\end{enumerate}
\end{theorem}

The next theorem makes explicit the relation between the
percolation full scaling limit and $SLE_6$.
Its proof (see Section~\ref{proofs}) relies on an inductive
procedure that makes use of $SLE_6$ at each step and allows
to obtain collections of loops with the correct distribution.

\begin{theorem} \label{main-result2}
A Continuum Nonsimple Loop process with the same distribution
as in Theorem~\ref{main-result1} can be constructed by 
a procedure in which each loop is obtained as the concatenation of an $SLE_6$
path with (a portion of) another $SLE_6$ path (see Figure~\ref{fig-sec3}).
This procedure is carried out first in a finite disk ${\mathbb D}_R$ of
radius $R$ in the plane (see Section~\ref{unit-disc}), and then an infinite
volume limit, ${\mathbb D}_R \to {\mathbb R}^2$, is taken.
\end{theorem}

\begin{remark} \label{conf-inv}
There are various possible ways to formulate the conformal
invariance properties of the Continuum Nonsimple Loop process.
One version is given in Theorem~\ref{thm-conf-inv}.
\end{remark}


Next we give some examples showing how the scaling limit of
various connectivity/crossing probabilities can be expressed
in terms of the loop process.
Although we cannot say whether this fact may eventually lead
to exact expressions going beyond Cardy's formula, it at least
shows that scaling limits of such probabilities exist and are
conformally invariant (early discussions of scaling limits of
connectivity functions and of the consequences of conformal
invariance for such quantities are given in~\cite{aizenman1,aizenman}).
The examples will also highlight the natural nested structure of the
collection of percolation cluster boundaries in the scaling limit.

\begin{figure}[!ht]
\begin{center}
\includegraphics[width=6cm]{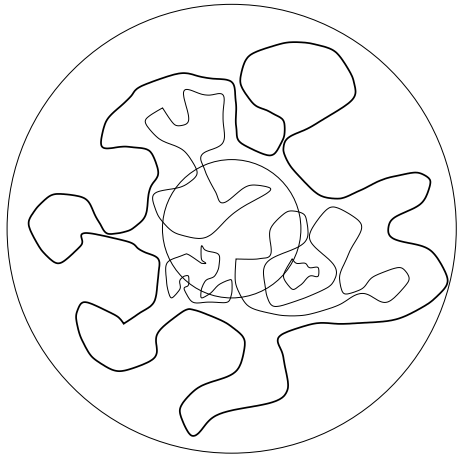}
\caption{An annulus whose inner disc is surrounded by a continuum
nonsimple loop.
There is no monochromatic crossing between the inner and outer discs.
Other continuum nonsimple loops are shown in the figure, but they do
not affect the connectivity between the inner and outer discs.}
\label{annulus}
\end{center}
\end{figure}

Consider first an annulus centered at $z$ with inner radius
$r_1$ and outer radius $r_2$ (see Figure~\ref{annulus}).
The scaling limit $p(r_1,r_2)$ of the probability of a crossing
of the annulus (by crossing here we refer to a ``monochromatic"
crossing, i.e., a crossing by either of the two colors, as discussed
in Section~\ref{lapa} -- see also Figure~\ref{loops}) can be expressed
as follows in terms of the random variable $N(z,r_1,r_2)$
defined in Theorem~\ref{features} above:
$p(r_1,r_2)$ is the probability
that $N(z,r_1,r_2)$ equals zero.
More generally, $N(z,r_1,r_2)$ represents the scaling limit of the
minimal number of cluster boundaries traversed by any path connecting
the inner and outer circles of the annulus.

\begin{figure}[!ht]
\begin{center}
\includegraphics[width=8cm]{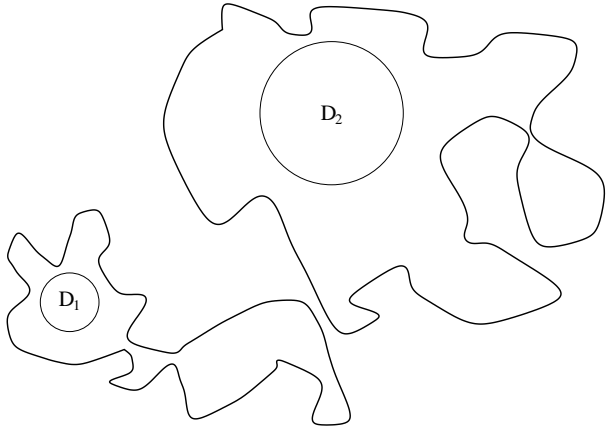}
\caption{Each one of the two disjoint discs in the figure is surrounded
by a continuum nonsimple loop that has the other disc in its exterior.
The minimal number of cluster boundaries that must be crossed to connect
the two discs is two.}
\label{discs}
\end{center}
\end{figure}

An example with more geometric structure involves two disjoint discs
$D_1$ and $D_2$ in the plane and the scaling limit $p(D_1,D_2)$ of
the probability that there is a crossing from $D_1$ to $D_2$ (see
Figure~\ref{discs}).
Here we let $N_1$ denote the number of distinct loops in the plane
that contain $D_1$ in their interior and $D_2$ in their exterior,
and define $N_2$ in the complementary way.
Then $p(D_1,D_2)$ is the probability that $N_1=N_2=0$, and the scaling
limit of the minimal number of cluster boundaries that must be crossed
to connect $D_1$ to $D_2$ is $N_1+N_2$.

One can also consider, as in~\cite{aizenman1,aizenman}, the probability
of a single monochromatic cluster in the exterior ${\cal E}$ of the union of
$m$ disjoint discs (or other regions) connecting all $m$ disc boundaries.
In the scaling limit, this can be expressed as the probability of the
event that there is a single continuous (nonsimple) curve in ${\cal E}$
touching all $m$ disc boundaries that {\it does not cross any of the loops
of the Continuum Nonsimple Loop process.}

We conclude this section by remarking that the Continuum Nonsimple Loop
process is just one example of a family of ``conformal loop ensembles"
that are related to $SLE$ and to the Gaussian Free Field
(see~\cite{werner3,werner5,werner6,ss,shw}), and are conjectured to describe
the full scaling limit of statistical mechanics models such as percolation,
Ising and Potts models.
The work of Lawler, Schramm, Sheffield and Werner has provided tools to define
such loop ensembles in the continuum and to study some of their properties.
But the key question concerning the proof of their connection to the discrete
models via a continuum scaling limit remains an important open challenge,
the only exception currently being percolation, as we show in this paper.

\section{Preliminary Definitions and Results} \label{defs}

We will find it convenient to identify the
real plane ${\mathbb R}^2$ and the complex plane $\mathbb C$.
We will also refer to the Riemann sphere ${\mathbb C} \cup \infty$
and the open upper half-plane $\mathbb H = \{ x+iy : y>0 \}$
(and its closure $\overline{\mathbb H}$),
where chordal $SLE$ will be defined (see Section~\ref{sle1}).
$\mathbb D$ will denote the open unit disc
${\mathbb D} = \{ z \in {\mathbb C} : |z|<1 \}$.

A domain $D$ of the complex plane $\mathbb C$ is a
nonempty, connected, open subset of $\mathbb C$; a
simply connected domain $D$ is said to be a Jordan
domain if its (topological) boundary $\partial D$
is a Jordan curve (i.e., a simple continuous loop).

We will make repeated use of Riemann's mapping theorem,
which states that if $D$ is any simply connected domain
other than the entire plane $\mathbb C$ and $z_0 \in D$,
then there is a unique conformal
map $f$ of $D$ onto $\mathbb D$ such that $f(z_0)=0$
and $f'(z_0)>0$.

\subsection{Compactification of ${\mathbb R}^2$}

When taking the scaling limit as the lattice mesh size
$\delta \to 0$ one can focus on fixed finite regions,
$\Lambda \subset {\mathbb R}^2$, or consider the whole
${\mathbb R}^2$ at once.
The second option avoids dealing with boundary conditions,
but requires an appropriate choice of metric.

A convenient way of dealing with the whole ${\mathbb R}^2$
is to replace the Euclidean metric with a distance function
$\Delta(\cdot,\cdot)$ defined on
${\mathbb R}^2 \times {\mathbb R}^2$ by
\begin{equation}
\Delta(u,v) = \inf_{\varphi} \int (1 + | {\varphi} |^2)^{-1} \, ds,
\end{equation}
where the infimum is over all smooth curves $\varphi(s)$
joining $u$ with $v$, parameterized by arclength $s$, and
where $|\cdot|$ denotes the Euclidean norm.
This metric is equivalent to the Euclidean metric in bounded
regions, but it has the advantage of making ${\mathbb R}^2$
precompact.
Adding a single point at infinity yields the compact space
$\dot{\mathbb R}^2$ which is isometric, via stereographic
projection, to the two-dimensional sphere.

\subsection{The Space of Curves} \label{space}

In dealing with the scaling limit we use the approach of
Aizenman-Burchard~\cite{ab}.
Denote by ${\cal S}_R$ the complete separable metric space
of continuous curves in the closure $\overline{\mathbb D}_R$
of the disc ${\mathbb D}_R$ of radius $R$ with the metric~(\ref{distance})
defined below.
Curves are regarded as equivalence classes of continuous
functions from the unit interval to $\overline{\mathbb D}_R$,
modulo monotonic reparametrizations.
$\gamma$ will represent a particular curve and $\gamma(t)$ a
parametrization of $\gamma$; ${\cal F}$ will represent a set
of curves (more precisely, a closed subset of ${\cal S}_R$).
$\text{d}(\cdot,\cdot)$ will denote the uniform metric
on curves, defined by
\begin{equation} \label{distance}
\text{d} (\gamma_1,\gamma_2) \equiv \inf
\sup_{t \in [0,1]} |\gamma_1(t) - \gamma_2(t)|,
\end{equation}
where the infimum is over all choices of parametrizations
of $\gamma_1$ and $\gamma_2$ from the interval $[0,1]$.
The distance between two closed sets of curves is defined
by the induced Hausdorff metric as follows:
\begin{equation} \label{hausdorff}
\text{dist}({\cal F},{\cal F}') \leq \varepsilon
\Leftrightarrow (\forall \, \gamma \in {\cal F}, \, \exists \,
\gamma' \in {\cal F}' \text{ with }
\text{d} (\gamma,\gamma') \leq \varepsilon,
\text{ and vice versa}).
\end{equation}
The space $\Omega_R$ of closed subsets of ${\cal S}_R$
(i.e., collections of curves in $\overline{\mathbb D}_R$)
with the metric~(\ref{hausdorff}) is also a complete
separable metric space.
We denote by ${\cal B}_R$ its Borel $\sigma$-algebra.

For each fixed $\delta>0$, the random curves that we consider
are polygonal paths
on the edges of the hexagonal lattice $\delta {\cal H}$,
dual to the triangular lattice $\delta {\cal T}$.
A superscript $\delta$ is added to indicate that the
curves correspond to a model with a ``short distance cutoff''
of magnitude $\delta$.

We will also consider the complete separable metric space ${\cal S}$
of continuous curves in $\dot{\mathbb R}^2$ with the distance
\begin{equation} \label{Distance}
\text{D} (\gamma_1,\gamma_2) \equiv \inf
\sup_{t \in [0,1]} \Delta(\gamma_1(t),\gamma_2(t)),
\end{equation}
where the infimum is again over all choices of parametrizations
of $\gamma_1$ and $\gamma_2$ from the interval $[0,1]$.
The distance between two closed sets of curves is again
defined by the induced Hausdorff metric as follows:
\begin{equation} \label{hausdorff-D}
\text{Dist}({\cal F},{\cal F}') \leq \varepsilon
\Leftrightarrow (\forall \, \gamma \in {\cal F}, \, \exists \,
\gamma' \in {\cal F}' \text{ with }
\text{D} (\gamma,\gamma') \leq \varepsilon
\text{ and vice versa}).
\end{equation}
The space $\Omega$ of closed sets of $\cal S$
(i.e., collections of curves in $\dot{\mathbb R}^2$)
with the metric~(\ref{hausdorff-D}) is also a complete
separable metric space.
We denote by ${\cal B}$ its Borel $\sigma$-algebra.

When we talk about convergence in distribution of random curves,
we always mean with respect to the uniform metric~(\ref{distance}),
while when we deal with closed collections of curves, we always
refer to the metric~(\ref{hausdorff}) or~(\ref{hausdorff-D}).

\begin{remark} \label{ab}
In this paper, the space $\Omega$ of closed sets of $\cal S$ is
used for collections of exploration paths (see Section~\ref{explo})
and cluster boundary loops and their scaling limits, $SLE_6$ paths
and continuum nonsimple loops.
\end{remark}

\subsection{Chordal $SLE$ in the Upper Half-Plane} \label{sle1}

The Stochastic Loewner Evolution ($SLE$) was introduced
by Schramm~\cite{schramm} as a tool for studying the
scaling limit of two-dimensional discrete (defined on a
lattice) probabilistic models whose scaling limits are
expected to be conformally invariant.
In this section we define the chordal version of $SLE$;
for more on the subject, the interested reader can consult
the original paper~\cite{schramm} as well as the fine
reviews by Lawler~\cite{lawler1}, Kager and Nienhuis~\cite{kn},
and Werner~\cite{werner4}, and Lawler's book~\cite{lawler2}.

Let $\mathbb H$ denote the upper half-plane.
For a given continuous real function $U_t$ with $U_0 = 0$,
define, for each $z \in \overline{\mathbb H}$, the function
$g_t(z)$ as the solution to the ODE
\begin{equation}
\partial_t g_t(z) = \frac{2}{g_t(z) - U_t},
\end{equation}
with $g_0(z) = z$.
This is well defined as long as $g_t(z) - U_t \neq 0$,
i.e., for all $t < T(z)$, where
\begin{equation}
T(z) \equiv \sup \{ t \geq 0 : \min_{s \in [0,t]} | g_s(z) - U_s| > 0 \}.
\end{equation}
Let $K_t \equiv \{ z \in \overline{\mathbb H} : T(z) \leq t \}$
and let ${\mathbb H}_t$ be the unbounded component of
${\mathbb H} \setminus K_t$; it can be shown that $K_t$ is bounded
and that $g_t$ is a conformal map from ${\mathbb H}_t$ onto $\mathbb H$.
For each $t$, it is possible to write $g_t(z)$ as
\begin{equation}
g_t(z) = z + \frac{2t}{z} + o(\frac{1}{z}),
\end{equation}
when $z \to \infty$.
The family $(K_t, t \geq 0)$ is called the {\bf Loewner chain}
associated to the driving function $(U_t, t \geq 0)$.

\begin{definition} \label{def-sle}
{\bf Chordal $SLE_{\kappa}$} is the Loewner chain $(K_t, t \geq 0)$
that is obtained when the driving function
$U_t = \sqrt{\kappa} B_t$ is $\sqrt{\kappa}$ times a standard
real-valued Brownian motion $(B_t, t \geq 0)$ with $B_0 = 0$.
\end{definition}

For all $\kappa \geq 0$, chordal $SLE_{\kappa}$ is almost surely generated
by a continuous random curve $\gamma$ in the sense that, for all $t \geq 0$,
${\mathbb H}_t \equiv {\mathbb H} \setminus K_t$ is the unbounded connected
component of ${\mathbb H} \setminus \gamma[0,t]$; $\gamma$ is called the
{\bf trace} of chordal $SLE_{\kappa}$.

\subsection{Chordal $SLE$ in 
a Jordan Domain} \label{sle2}
Let $D \subset {\mathbb C}$ 
be a Jordan domain.
By Riemann's mapping theorem, there are (many) conformal maps
from the upper half-plane $\mathbb H$ onto $D$.
In particular, given two distinct points $a,b \in \partial D$,
there exists a
conformal map $f$ from $\mathbb H$ onto $D$ such that $f(0)=a$
and $f(\infty) \equiv \lim_{|z| \to \infty} f(z) = b$.
In fact, the choice of the points $a$ and $b$ on the boundary
of $D$ only characterizes $f(\cdot)$ up to a multiplicative
factor, since $f(\lambda \, \cdot)$ would also do.

Suppose that $(K_t, t \geq 0)$ is a chordal $SLE_{\kappa}$ in
$\mathbb H$ as defined above; we define chordal $SLE_{\kappa}$
$(\tilde K_t, t \geq 0)$ in $D$ from $a$ to $b$ as the
image of the Loewner chain $(K_t, t \geq 0)$ under $f$.
It is possible to show, using scaling properties of
$SLE_{\kappa}$, that the law of $(\tilde K_t, t \geq 0)$
is unchanged, up to a linear time-change, if we replace
$f(\cdot)$ by $f(\lambda \, \cdot)$.
This makes it natural to consider $(\tilde K_t, t \geq 0)$ as
a process from $a$ to $b$ in $D$, ignoring the role of $f$.

We are interested in the case $\kappa = 6$, for which
$(K_t, t \geq 0)$ is generated by a continuous, nonsimple,
non-self-crossing curve $\gamma$ with Hausdorff dimension $7/4$.
We will denote by $\gamma_{D,a,b}$ the image of $\gamma$ under $f$
and call it the trace of chordal $SLE_6$ in $D$ from $a$ to $b$;
$\gamma_{D,a,b}$ is a continuous nonsimple curve inside $D$ from $a$
to $b$, and it can be given a parametrization $\gamma_{D,a,b}(t)$
such that $\gamma_{D,a,b}(0)=a$ and $\gamma_{D,a,b}(1)=b$, so that
we are in the metric framework described in Section~\ref{space}.
It will be convenient to think of $\gamma_{D,a,b}$ as an
oriented path, with orientation from $a$ to $b$.

\subsection{Rad\'o's Theorem} \label{rado}

We present here Rad\'o's theorem~\cite{rado} (see also Theorem~2.11
of~\cite{pommerenke}), which deals with sequences of Jordan domains
and the corresponding conformal maps from the unit disc, and will
be used in the proof of the key Lemma~\ref{strong-smirnov}.

Since the theorem deals with Jordan domains, the conformal maps
from the unit disc to those domains have a continuous extension
to ${\mathbb D} \cup \partial {\mathbb D}$.
With a slight abuse of notation, we do not distinguish between
the conformal maps and their continuous extensions.
\begin{theorem} \label{rado-thm}
For $k=1,2,\ldots$ , let $J_k$ and $J$ be Jordan curves parameterized
respectively by $\phi_k(t)$ and $\phi(t)$, $t \in [0,1]$, and let
$f_k$ and $f$ be conformal maps from $\mathbb D$ onto the inner
domains of $J_k$ and $J$ such that $f_k(0)=f(0)$ and $f_k'(0)>0$,
$f'(0)>0$ for all $k$.
If $\phi_k \to \phi$ as $k \to \infty$ uniformly in $[0,1]$ then
$f_k \to f$ as $k \to \infty$ uniformly in $\overline{\mathbb D}$.
\end{theorem}

The type of convergence of sequences of Jordan domains $\{ D_k \}$
to a Jordan domain $D$ that will be encountered in Lemma~\ref{strong-smirnov}
is such that $\partial D_k$ converges, as $k \to \infty$, to
$\partial D$ in the uniform metric~(\ref{distance}) on continuous
curves, which is clearly sufficient to apply Theorem~\ref{rado-thm}.

\section{Construction of the Continuum Nonsimple Loops}
\label{construct}

\subsection{Construction of a Single Loop} \label{single-loop}

\begin{figure}[!ht]
\begin{center}
\includegraphics[width=8cm]{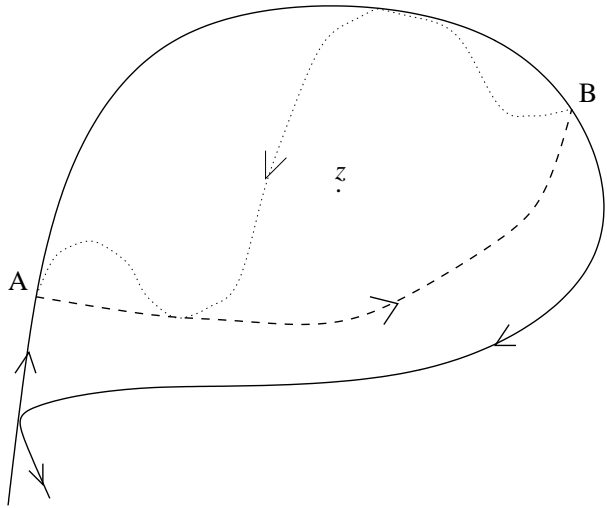}
\caption{Construction of a continuum loop around $z$ in three steps.
A domain $D$ is formed by the solid curve. The dashed curve is an
excursion $\cal E$ (from A to B) of an $SLE_6$ in $D$ that creates
a subdomain $D'$ containing $z$. The dotted curve $\gamma'$ is an
$SLE_6$ in $D'$ from B to A. A loop is formed by $\cal E$ followed
by $\gamma'$.}
\label{fig-sec3}
\end{center}
\end{figure}

As a preview to the full construction, we explain how to
construct a single loop using two $SLE_6$ paths inside a
domain $D$ whose boundary is assumed to have a given
orientation (clockwise or counterclockwise).
This is done in three steps (see Figure~\ref{fig-sec3}), of
which the first consists in choosing two points $a$ and $b$
on the boundary $\partial D$ of $D$ and ``running'' a chordal
$SLE_6$, $\gamma = \gamma_{D,a,b}$, from $a$ to $b$ inside $D$.
As explained in Section~\ref{sle2}, we consider $\gamma$
as an oriented path, with orientation from $a$ to $b$.
The set $D \setminus \gamma_{D,a,b}[0,1]$ is a countable
union of its connected components, which are open and simply
connected.
If $z$ is a deterministic point in $D$, then with probability
one, $z$ is not touched by $\gamma$~\cite{rs} and so it belongs
to a unique domain in $D \setminus \gamma_{D,a,b}[0,1]$
that we denote $D_{a,b}(z)$.

\begin{figure}[!ht]
\begin{center}
\includegraphics[width=6cm]{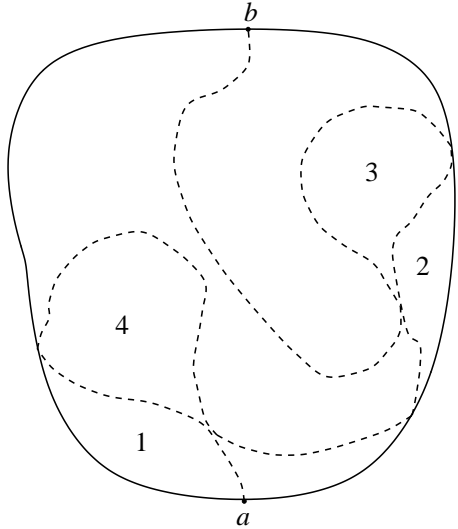}
\caption{Schematic diagram showing the four types of (sub)domains
formed by a dashed curve $\gamma$ from $a$ to $b$ inside a domain
whose boundary is the solid curve.}
\label{four-types}
\end{center}
\end{figure}

The elements of $D \setminus \gamma_{D,a,b}[0,1]$ can
be conveniently characterized in terms of how a point $z$ in
the interior of the component was first ``trapped'' at some
time $t_1$ by $\gamma[0,t_1]$, perhaps together with either
$\partial_{a,b} D$ or $\partial_{b,a} D$
(the portions of the boundary $\partial D$ from $a$ to $b$
counterclockwise or clockwise respectively) --- see Figure~\ref{four-types}:
(1) those components whose boundary contains a segment of
$\partial_{b,a} D$ between two successive visits at
$\gamma_0(z)=\gamma(t_0)$ and $\gamma_1(z)=\gamma(t_1)$ to
$\partial_{b,a} D$ (where here and below $t_0<t_1$), (2) the
analogous components with $\partial_{b,a} D$  replaced by the
other part of the boundary $\partial_{a,b} D$, (3) those components
formed when $\gamma_0(z)=\gamma(t_0)=\gamma(t_1)=\gamma_1(z) \in D$
with $\gamma$ winding about $z$ in a counterclockwise direction
between $t_0$ and $t_1$, and finally (4) the analogous clockwise
components.

We give to the boundary of a domain of type~3 or 4 the orientation
induced by how the curve $\gamma$ winds around the points inside
that domain.
For a domain $D' \ni z$ of type~1 or 2 which is produced by an
``excursion" $\cal E$ from $\gamma_0(z) \in \partial D$ to
$\gamma_1(z) \in \partial D$, the part of the boundary that
corresponds to the inner perimeter of the excursion $\cal E$
(i.e., the perimeter of $\gamma$ seen from $z$) is oriented
according to the direction of $\gamma$, i.e., from $\gamma_0(z)$
to $\gamma_1(z)$.

If we assume that $\partial D$ is oriented from $a$ to $b$
clockwise, then the boundaries of domains of type 2 have
a well defined orientation, while the boundaries of domains of
type 1 do not, since they are composed of two parts which are
both oriented from the beginning to the end of the excursion
that produced the domain.

Now, let $D'$ be a domain of type 1 and let $A$ and $B$ be
respectively the starting and ending point of the excursion
that generated $D'$.
The second step to construct a loop is to run a chordal
$SLE_6$, $\gamma' = \gamma_{D',B,A}$, inside $D'$ from
$B$ to $A$; the third and final step consists in pasting
together $\cal E$ and $\gamma'$.

Running $\gamma'$ inside $D'$ from $B$ to $A$
partitions $D' \setminus \gamma'$ into new domains.
Notice that if we assign an orientation to the boundaries of
these domains according to the same rules used above, all of
those boundaries have a well defined orientation, so that
the construction of loops just presented can be iterated
inside each one of these domains (as well as inside each
of the domains of type~2, 3 and 4 generated by $\gamma_{D,a,b}$
in the first step).
This will be done in the next section.

%
%


\subsection{The Full Construction Inside The Unit Disc}
\label{unit-disc}

In this section we define the Continuum Nonsimple Loop
process inside the unit disc  ${\mathbb D} = {\mathbb D}_1$
via an inductive procedure.
Later, in order to define the continuum nonsimple loops in
the whole plane, the unit disc will be replaced by a growing
sequence of large discs, ${\mathbb D}_R$, with $R\to \infty$
(see Theorem~\ref{thm-therm-lim}).
The basic ingredient in the algorithmic construction, given
in the previous section, consists of a chordal $SLE_6$ path
$\gamma_{D,a,b}$ between two points $a$ and $b$ of the boundary
$\partial D$ of a given simply connected domain $D \subset {\mathbb C}$.

We will organize the inductive procedure in steps,
each one corresponding to one $SLE_6$ inside a certain domain
generated by the previous steps.
To do that, we need to order the domains present at the end
of each step, so as to choose the one to use in the next step.
For this purpose, we introduce a deterministic countable
set of points $\cal P$ that are dense in $\mathbb C$ and
are endowed with a deterministic order (here and below
by deterministic we mean that they are assigned before
the beginning of the construction and are independent
of the $SLE_6$'s).

The first step consists of an $SLE_6$ path,
$\gamma_1 = \gamma_{{\mathbb D},-i,i}$, inside $\mathbb D$
from $-i$ to $i$, which produces many domains that are
the  connected components of the set
$\mathbb D \setminus \gamma_1 [0,1]$.
These domains can be priority-ordered according to the
maximal $x$- or $y$- coordinate distances between points on
their boundaries and using the rank of the points in $\cal P$
(contained in the domains) to break ties, as follows.
For a domain $D$, let $\text{d}_m(D)$ be the maximal $x$-
or $y$-distance between points on its boundary, whichever
is greater.
Domains with larger $\text{d}_m$ have higher priority, and
if two domains have the same $\text{d}_m$, the one containing
the highest ranking point of $\cal P$ from those two domains
has higher priority.
The priority order of domains of course changes as the
construction proceeds and new domains are formed.

The second step of the construction consists of an $SLE_6$
path, $\gamma_2$, that is produced in the domain with highest
priority (after the first step).
Since all the domains that are produced in the construction
are Jordan domains, as explained in the discussion following
Corollary~\ref{jordan}, for all steps we can use the definition
of chordal $SLE$ given in Section~\ref{sle2}.

As a result of the construction, the $SLE_6$ paths are
naturally ordered: $\{ \gamma_j \}_{j \in {\mathbb N}}$.
It will be shown (see especially the proof of Theorem~\ref{thm-convergence}
below) that every domain that is formed during the construction
is eventually used (this is in fact one important requirement in
deciding how to order the domains and therefore how to organize
the construction).

So far we have not explained how to choose the starting
and ending points of the $SLE_6$ paths on the boundaries
of the domains.
In order to do this, we give an orientation to the
boundaries of the domains produced by the construction
according to the rules explained in Section~\ref{single-loop}.
We call {\bf monochromatic} a boundary which gets, as
a consequence of those rules, a well defined (clockwise
or counterclockwise) orientation; the choice of this term
will be clarified when we discuss the lattice version of
the loop construction below.
We will generally take our initial domain ${\mathbb D}_1$
(or ${\mathbb D}_R$) to have a monochromatic boundary
(either clockwise or counterclockwise orientation).

It is easy to see by induction that the boundaries that
are not monochromatic are composed of two ``pieces'' joined
at two special points (call them A and B, as in the example
of Section~\ref{single-loop}), such that one piece is a
portion of the boundary of a previous domain, and the
other is the inner perimeter of an excursion (see again
Section~\ref{single-loop}).
Both pieces are oriented in the same direction, say from A
to B (see Figure~\ref{fig-sec3}).

For a domain whose boundary is not monochromatic, we make
the ``natural'' choice of starting and ending points,
corresponding to the end and beginning of the excursion
that produced the domain (the points B and A respectively,
in the example above).
As explained in Section~\ref{single-loop}, when such
a domain is used with this choice of points on the
boundary, a loop is produced, together with other domains,
whose boundaries are all monochromatic.

For a domain whose boundary is monochromatic, and therefore
has a well defined orientation, there are various procedures
which would yield the ``correct'' distribution for the resulting
Continuum Nonsimple Loop process; one possibility is as follows.

Given a domain $D$, $a$ and $b$ are chosen so that, of all pairs
$(u,v)$ of points in $\partial D$, they maximize $|\text{Re}(u-v)|$
if $|\text{Re}(u-v)| \geq |\text{Im}(u-v)|$, or else they maximize
$|\text{Im}(u-v)|$.
If the choice is not unique, to restrict the number of pairs
one looks at those pairs, among the ones already obtained, that
maximize the other of $\{ |\text{Re}(u-v)|, |\text{Im}(u-v)| \}$.
Notice that this leaves at most two pairs of points; if that's
the case, the pair that contains the point with minimal real
(and, if necessary, imaginary) part is chosen.
The iterative procedure produces a loop every time a domain
whose boundary is not monochromatic is used.
Our basic loop process consists of the collection of all loops
generated by this inductive procedure (i.e., the limiting object
obtained from the construction by letting the number of steps
$k \to \infty$), to which we add a ``trivial" loop for each $z$
in $\mathbb D$, so  that the collection of loops is closed in
the appropriate sense~\cite{ab}.
The Continuum Nonsimple Loop process in the whole plane is
introduced in Theorem~\ref{thm-therm-lim}, Section~\ref{main-results}.
There, a ``trivial" loop for each $z \in {\mathbb C} \cup \infty$
has to be added to make the space of loops closed.

\section{Lattices and Paths} \label{lapa}

We will denote by $\cal T$ the two-dimensional triangular lattice,
whose sites we think of as the elementary cells of a regular hexagonal
lattice $\cal H$ embedded in the plane as in Figure~\ref{loops}.
Two hexagons are {\bf neighbors} if they are adjacent, i.e., if they
have a common edge.
A sequence $(\xi_0, \ldots, \xi_n)$ of hexagons 
such that $\xi_{i-1}$ and $\xi_i$ are neighbors 
for all $i= 1, \ldots, n$ and $\xi_i \neq \xi_j$ whenever $i \neq j$
will be called a {\bf $\cal T$-path} and denoted by $\pi$.
If the first and last sites of the path are neighbors, 
the path will be called a {\bf $\cal T$-loop}.

A finite set $D$ of hexagons is {\bf connected} if any two hexagons
in $D$ can be joined by a $\cal T$-path contained in $D$.
We say that a finite set $D$ of hexagons 
is {\bf simply connected} if both $D$ and its complement 
are connected. 
For a simply connected set $D$ of hexagons, we denote by
$\Delta D$ its {\bf external site boundary}, or {\bf s-boundary}
(i.e., the set of hexagons that do not belong to $D$
but are adjacent to hexagons in $D$), and by $\partial D$ the
topological boundary of $D$ when $D$ is considered as a domain of
$\mathbb C$.
We will call a bounded, simply connected subset $D$ of $\cal T$
a {\bf Jordan set} if its s-boundary $\Delta D$ is a $\cal T$-loop.

For a Jordan set $D \subset {\cal T}$, a vertex $x \in {\cal H}$
that belongs to $\partial D$ can be either of two types, according to whether
the edge incident on $x$ that is not in $\partial D$ belongs to a hexagon
in $D$ or not.
We call a vertex of the second type an {\bf e-vertex} (e for ``external''
or ``exposed'').

Given a Jordan set $D$ and two e-vertices $x,y$ in $\partial D$,
we denote by $\partial_{x,y} D$ the portion of $\partial D$
traversed counterclockwise from $x$ to $y$, and call it the
{\bf right boundary}; the remaining part of the boundary is
denote by $\partial_{y,x} D$ and is called the {\bf left boundary}.
Analogously, the portion of $\Delta_{x,y} D$ of $\Delta D$ whose
hexagons are adjacent to $\partial_{x,y} D$ is called the
{\bf right s-boundary} and the remaining part the {\bf left s-boundary}.


A {\bf percolation configuration}
$\sigma = \{ \sigma(\xi) \}_{\xi \in \cal T} \in \{ -1, +1 \}^{\cal T}$
on $\cal T$ is an assignment of $-1$ (equivalently, yellow) or $+1$
(blue) to each site of $\cal T$ (i.e., to each hexagon of
$\cal H$ -- see Figure~\ref{loops}). 
For a domain $D$ of the plane, the restriction to the subset
$D \cap \cal T$ of $\cal T$ of the percolation configuration
$\sigma$ is denoted by $\sigma_D$.
On the space of configurations $\Sigma = \{ -1,+1 \}^{\cal T}$,
we consider the usual product topology and denote by $\mathbb P$
the uniform measure, corresponding to Bernoulli percolation with
equal density of yellow (minus) and blue (plus) hexagons, which is
critical percolation in the case of the triangular lattice.

A (percolation) {\bf cluster} is a maximal, connected, monochromatic
subset of $\cal T$; we will distinguish between blue (plus) and
yellow (minus) clusters.
The {\bf boundary} of a cluster $D$ is the set of edges of
$\cal H$ that surround the cluster (i.e., its Peierls contour);
it coincides with the topological boundary of $D$ considered as a
domain of $\mathbb C$.
The set of all boundaries is a collection of ``nested'' simple loops
along the edges of $\cal H$.

Given a percolation configuration $\sigma$, we associate an
arrow to each edge of $\cal H$ belonging to the boundary of
a cluster in such a way that the hexagon to the right of the edge
with respect to the direction of the arrow is blue (plus).
The set of all boundaries then becomes a collection of nested,
oriented, simple loops.
A {\bf boundary path} (or {\bf b-path}) $\gamma$ is a sequence
$(e_0, \ldots, e_n)$ of distinct edges of $\cal H$ belonging
to the boundary of a cluster and such that $e_{i-1}$ and $e_i$
meet at a vertex of $\cal H$ for all $i= 1, \ldots, n$.
To each b-path, we can associate a direction according to the
direction of the edges in the path.

Given a b-path $\gamma$, we denote by $\Gamma_B(\gamma)$
(respectively, $\Gamma_Y(\gamma)$) the set of blue (resp.,
yellow) hexagons (i.e., sites of $\cal T$) adjacent to
$\gamma$; we also let
$\Gamma(\gamma) \equiv \Gamma_B(\gamma) \cup \Gamma_Y(\gamma)$.

\subsection{The Percolation Exploration Process and Path} \label{explo}

\begin{figure}[!ht]
\begin{center}
\includegraphics[width=6cm]{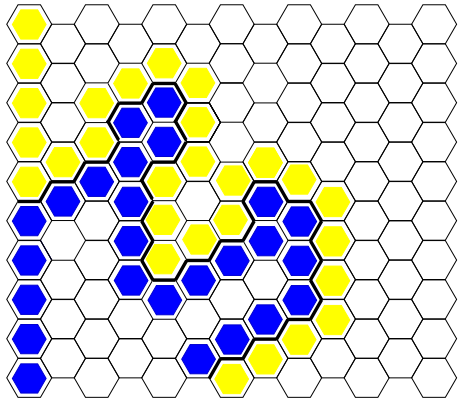}
\caption{Percolation exploration process in a portion of the
hexagonal lattice with $\pm$ boundary conditions on the first
column, corresponding to the boundary of the region where the
exploration is carried out.
The colored hexagons that do not belong to the first column
have been ``explored'' during the exploration process.
The heavy line between yellow (light) and blue (dark) hexagons
is the exploration path produced by the exploration process.}
\label{fig2-sec4}
\end{center}
\end{figure}

For a Jordan set $D \subset {\cal T}$ and two
e-vertices $x,y$ in $\partial D$, imagine coloring blue
all the hexagons in $\Delta_{x,y} D$ and yellow all those
in $\Delta_{y,x} D$.
Then, for any percolation configuration $\sigma_D$ inside
$D$, there is a unique b-path $\gamma$ from $x$ to $y$ which
separates the blue cluster adjacent to $\Delta_{x,y} D$ from
the yellow cluster adjacent to $\Delta_{y,x} D$.
We call $\gamma = \gamma_{D,x,y}(\sigma_D)$ a
{\bf percolation exploration path} (see Figure~\ref{fig2-sec4}).

An exploration path $\gamma$ can be decomposed into
{\bf left excursions} $\cal E$, i.e., maximal b-subpaths
of $\gamma$ that do not use edges of the left boundary $\partial_{y,x} D$.
Successive left excursions are separated by portions of
$\gamma$ that contain only edges of the left boundary $\partial_{y,x} D$.
Analogously, $\gamma$ can be decomposed into {\bf right excursions},
i.e., maximal b-subpaths of $\gamma$ that do not use edges of the right
boundary $\partial_{x,y} D$.
Successive right excursions are separated by portions of
$\gamma$ that contain only edges of the right boundary $\partial_{x,y} D$.

Notice that the exploration path
$\gamma=\gamma_{D,x,y}(\sigma_D)$ only depends on
the percolation configuration $\sigma_D$ inside $D$
and the positions of the e-vertices $x$ and $y$;
in particular, it does not depend on the color of
the hexagons in $\Delta D$, since it is defined by
imposing fictitious $\pm$ boundary conditions on $D$.
To see this more clearly, we next show how to construct
the percolation exploration path dynamically, via
the {\bf percolation exploration process} defined
below.

Given a Jordan set $D \subset {\cal T}$ and two
e-vertices $x,y$ in $\partial D$, assign to $\partial_{x,y} D$
a counterclockwise orientation (i.e., from $x$ to $y$) and
to $\partial_{y,x} D$ a clockwise orientation.
Call $e_x$ the edge incident on $x$ that does not belong to
$\partial D$ and orient it in the direction of $x$;
this is the ``starting edge'' of an exploration procedure
that will produce an oriented path inside $D$ along
the edges of $\cal H$, together with two \emph{nonsimple}
monochromatic paths on $\cal T$.
>From $e_x$, the process moves along the edges of
hexagons in $D$ according to the rules below.
At each step there are two possible edges (left or right edge
with respect to the current direction of exploration) to choose
from, both belonging to the same hexagon $\xi$ contained in $D$
or $\Delta D$.

\begin{itemize}
\item If $\xi$ belongs to $D$ and has not been previously ``explored,"
its color is determined by flipping a fair coin and then
the edge to the left (with respect to the direction in which
the exploration is moving) is chosen if $\xi$ is blue (plus),
or the edge to the right is chosen if $\xi$ is yellow (minus).
\item If $\xi$ belongs to $D$ and has been previously explored, the
color already assigned to it is used to choose an edge according to
the rule above.
\item If $\xi$ belongs to the right external boundary $\Delta_{x,y} D$,
the left edge is chosen.
\item If $\xi$ belongs to the left external boundary $\Delta_{y,x} D$,
the right edge is chosen.
\item The exploration process stops when it reaches $b$.
\end{itemize}

We can assign an arrow to each edge in the path in such
a way that the hexagon to the right of the edge with respect
to the arrow is blue; for edges in $\partial D$, we assign the
arrows according to the direction assigned to the boundary.
In this way, we get an oriented path, whose shape and orientation
depend solely on the color of the hexagons explored during
the construction of the path.

%


When we present the discrete construction, we will encounter
Jordan sets $D$ with two e-vertices $x,y \in \partial D$
assigned in some way to be discussed later.
Such domains will have either monochromatic (plus or
minus) boundaries or $\pm$ boundary conditions,
corresponding to having both $\Delta_{x,y} D$ and
$\Delta_{y,x} D$ monochromatic, but of different
colors.

As explained, the exploration path $\gamma_{D,x,y}$
does not depend on the color of $\Delta D$, but the
interpretation of $\gamma_{D,x,y}$ does.
For domains with $\pm$ boundary conditions, the
exploration path represents the interface between the
yellow cluster containing the yellow portion of the
s-boundary of $D$ and the blue cluster containing its
blue portion.

For domains with monochromatic blue (resp., yellow)
boundary conditions, the exploration path represents portions
of the boundaries of yellow (resp., blue)
clusters touching $\partial_{y,x} D$ and adjacent to
blue (resp., yellow) hexagons that are the starting point
of a blue (resp., yellow) path (possibly an empty path)
that reaches $\partial_{x,y} D$, pasted together using
portions of $\partial_{y,x} D$.

In order to study the continuum scaling limit of an exploration
path, we introduce the following definitions.
\begin{definition} \label{approx}
Given a Jordan domain $D$ of the plane,
we denote by $D^{\delta}$ the largest Jordan set of hexagons
of the scaled hexagonal lattice $\delta {\cal H}$ that is
contained in $D$, and call it the {\bf $\delta$-approximation}
of $D$.
\end{definition}

It is clear that
$\partial D^{\delta}$ converges to $\partial D$ in the
metric~(\ref{distance}).

\begin{definition} \label{exp-path}
Let $D$ be a Jordan domain of the plane and $D^{\delta}$
its $\delta$-approximation.
For $a,b \in \partial D$, choose the pair $(x_a,x_b)$
of e-vertices in $\partial D^{\delta}$ closest to, respectively,
$a$ and $b$ (if there are two such vertices closest to $a$,
we choose, say, the first one encountered going clockwise
along $\partial D^{\delta}$, and analogously for $b$).
Given a percolation configuration $\sigma$, we define
the {\bf exploration path} $\gamma^{\delta}_{D,a,b}(\sigma)
\equiv \gamma_{D^{\delta},x_a,x_b}(\sigma)$.
\end{definition}

For a fixed $\delta>0$, the measure $\mathbb P$ on percolation
configurations $\sigma$ induces a measure $\mu^{\delta}_{D,a,b}$
on exploration paths $\gamma^{\delta}_{D,a,b}(\sigma)$.
In the continuum scaling limit, $\delta \to 0$, one is
interested in the weak convergence of $\mu^{\delta}_{D,a,b}$ to a
measure $\mu_{D,a,b}$ supported on continuous curves, with respect
to the uniform metric~(\ref{distance}) on continuous curves.

One of the main tools in this paper is the result on convergence
to $SLE_6$ announced by Smirnov~\cite{smirnov} (see also~\cite{smirnov-long}),
whose detailed proof is to appear~\cite{smirnov1}: \emph{The distribution
of $\gamma^{\delta}_{D,a,b}$ converges, as $\delta \to 0$, to that
of the trace of chordal $SLE_6$ inside $D$ from $a$ to $b$, with
respect to the uniform metric~(\ref{distance}) on continuous curves}.

Actually, we will rather use a slightly stronger conclusion,
given as statement (S) at the beginning of Section~\ref{main-results}
below, a version of which, according to~\cite{sw} (see p.~734 there),
and~\cite{smirnov-private-comm}, will be contained in~\cite{smirnov1}.
This stronger statement is that the convergence of the percolation process
to $SLE_6$ takes place \emph{locally uniformly} with respect to the shape
of the domain $D$ and the positions of the starting and ending points $a$
and $b$ on its boundary $\partial D$.
We will use this version of convergence to $SLE_6$ to identify the
Continuum Nonsimple Loop process with the scaling limit of \emph{all}
critical percolation clusters.
A detailed proof of statement (S) can be found
in~\cite{cn2}.
Although the convergence statement in (S) is stronger than those
in~\cite{smirnov,smirnov-long}, we note that it is restricted to
Jordan domains, a restriction not present in~\cite{smirnov,smirnov-long}.



Before concluding this section, we give one more definition.
Consider the exploration path $\gamma = \gamma^{\delta}_{D,x,y}$
and the set $\Gamma(\gamma) = \Gamma_Y(\gamma) \cup \Gamma_B(\gamma)$.
The set $D^{\delta} \setminus \Gamma(\gamma)$ is the
union of its connected components (in the lattice sense),
which are simply connected.
If the domain $D$ is large and the e-vertices
$x_a, y_a \in \partial D^{\delta}$ are not too close
to each other, then with high probability the exploration
process inside $D^{\delta}$ will make large excursions
into $D^{\delta}$, so that
$D^{\delta} \setminus \Gamma(\gamma)$ will have more
than one component.
Given a point $z \in {\mathbb C}$ contained in
$D^{\delta} \setminus \Gamma(\gamma)$,
we will denote by $D^{\delta}_{a,b}(z)$ the domain
corresponding to the unique element of
$D^{\delta} \setminus \Gamma(\gamma)$ that contains
$z$ (notice that for a deterministic
$z \in D$, $D^{\delta}_{a,b}(z)$ is well defined
with high probability for $\delta$ small, i.e., when
$z \in D^{\delta}$ and $z \notin \Gamma(\gamma)$).

\subsection{Discrete Loop Construction}

\begin{figure}[!ht]
\begin{center}
\includegraphics[width=8cm]{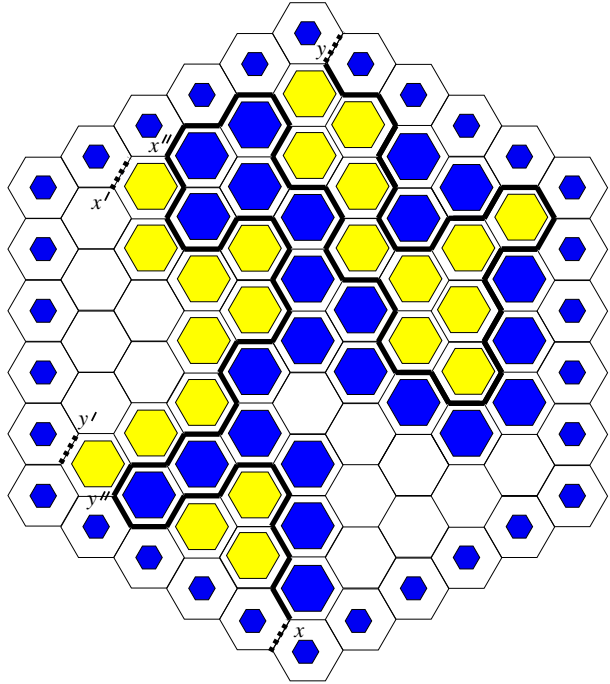}
\caption{First step of the construction of the outer contour
of a cluster of yellow/minus (light in the figure) hexagons
consisting of an exploration (heavy line) from the e-vertex $x$
to the e-vertex $y$.
The ``starting edge" and ``ending edge" of the exploration
path are indicated by dotted segments next to $x$ and $y$.
The outer layer of hexagons does not belong to the domain
where the explorations are carried out, but represents its
monochromatic blue/plus external site boundary.
$x''$ and $y''$ are the ending and starting points of a left
excursion that determines a new domain $D'$, and $x'$ and $y'$
are the vertices where the edges that separate the yellow and
blue portions of the s-boundary of $D'$ intersect $\partial D'$.
$x'$ and $y'$ will be respectively the beginning and end of a
new exploration path whose ``starting edge" and ``ending edge"
are indicated by dotted segments next to those points.}
\label{fig3-sec4}
\end{center}
\end{figure}

\begin{figure}[!ht]
\begin{center}
\includegraphics[width=8cm]{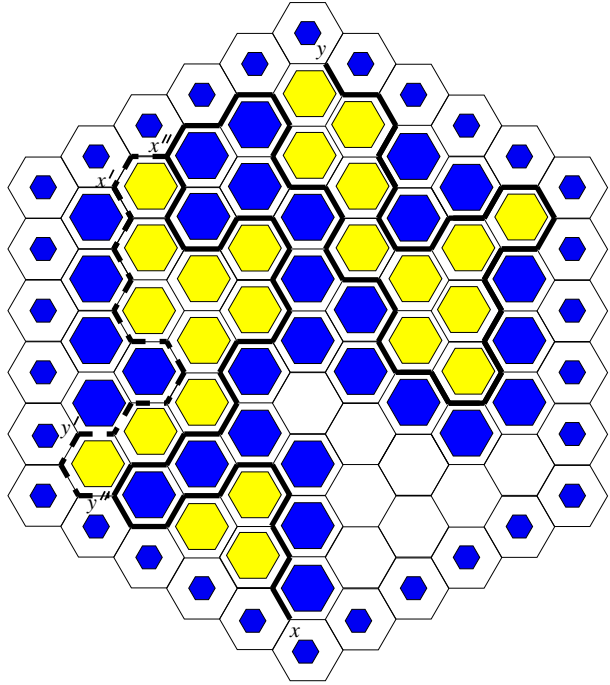}
\caption{Second step of the construction of the outer contour
of a cluster of yellow/minus (light in the figure) hexagons
consisting of an exploration from $x'$ to $y'$ whose resulting
path (heavy broken line) is pasted to the left excursion generated
by the previous exploration with the help of edges (indicated again
by a heavy broken line) of $\partial D$ contained between $x'$ and
$x''$ and between $y'$ and $y''$.}
\label{fig4-sec4}
\end{center}
\end{figure}

Next, we show how to construct, by twice using the
exploration process described in Section~\ref{explo},
a loop $\Lambda$ along the edges of ${\cal H}$
corresponding to the external boundary of a monochromatic
cluster contained in a large, simply connected, Jordan
set $D$ with monochromatic blue (say) boundary conditions
(see Figures~\ref{fig3-sec4} and~\ref{fig4-sec4}).

Consider the exploration path $\gamma = \gamma_{D,x,y}$
and the sets $\Gamma_Y(\gamma)$ and $\Gamma_B(\gamma)$
(see Figure~\ref{fig3-sec4}).
The set $D \setminus \{ \Gamma_Y(\gamma) \cup  \Gamma_B(\gamma) \}$
is the union of its connected components (in the lattice sense),
which are simply connected.
If the domain $D$ is large and the e-vertices
$x, y \in \partial D$ are chosen not too close
to each other, with large probability the exploration
process inside $D$ will make large excursions into
$D$, so that
$D \setminus \{ \Gamma_Y(\gamma) \cup  \Gamma_B(\gamma) \}$
will have many components.

There are four types of components which may be usefully
thought of in terms of their external site boundaries:
(1) those components  whose  site boundary contains both
sites in $\Gamma_Y(\gamma)$ and  $\Delta_{y,x} D$, (2)
the analogous components  with $\Delta_{y,x} D$ replaced
by  $\Delta_{x,y} D$ and $\Gamma_Y(\gamma)$ by
$\Gamma_B(\gamma)$, (3) those components whose site
boundary only contains sites  in $\Gamma_Y(\delta)$,
and finally (4) the analogous components  with
$\Gamma_Y(\gamma)$  replaced by $\Gamma_B(\gamma)$.

Notice that the components of type~1 are the only
ones with $\pm$ boundary conditions, while all other
components have monochromatic s-boundaries.
For a given component $D'$ of type~1, we can identify
the two edges that separate the yellow and blue portions
of its s-boundary.
The vertices $x'$ and $y'$ of $\cal H$ where those
two edges intersect $\partial D'$ are e-vertices and
are chosen to be the starting and ending points of
the exploration path $\gamma_{D',x',y'}$ inside $D'$.

If $x'', y'' \in \partial D$ are respectively the ending
and starting points of the left excursion $\cal E$ of $\gamma_{D,x,y}$
that ``created" $D'$, by pasting together $\cal E$ and
$\gamma_{D',x',y'}$ with the help of the edges of $\partial D$
contained between $x'$ and $x''$ and between $y'$ and $y''$,
we get a loop $\Lambda$ which corresponds to the
boundary of a yellow cluster adjacent to $\partial_{y,x} D$
(see Figure~\ref{fig4-sec4}).
Notice that the path $\gamma_{D',x',y'}$ in general splits
$D'$ into various other domains, all of which have monochromatic
boundary conditions.

%
%

\subsection{Full Discrete Construction} \label{full}

We now give the algorithmic construction for discrete
percolation  which is the analogue of the continuum one.
Each step of the construction is a single percolation
exploration  process; the order of successive steps is
organized as in the continuum construction detailed in
Section~\ref{unit-disc}.
We start with the largest Jordan set
$D^{\delta}_0 = {\mathbb D}^{\delta}$
of hexagons that is contained in the unit disc $\mathbb D$.
We will also make use of the countable set $\cal P$ of
points dense in $\mathbb C$ that was introduced earlier.

The first step consists of an exploration process inside
$D^{\delta}_0$.
For this, we need to select two points $x$ and $y$ in
$\partial D^{\delta}_0$ (which identify the starting and
ending edges).
We choose for $x$ the e-vertex closest to $-i$, and for
$y$ the e-vertex closest to $i$ (if there are two such
vertices closest to $-i$, we can choose, say, the one with
smallest real part, and analogously for $i$).
The first exploration produces a path $\gamma^{\delta}_1$
and, for $\delta$ small, many new domains of all four types.
These domains are ordered according to the maximal $x$- or
$y$- distance $\text{d}_m$ between points on their boundaries
and, if necessary, with the help of points in $\cal P$, as in
the continuum case, and that order is used, at each step of the
construction, to determine the next exploration process.
With this choice, the exploration processes and paths are
naturally ordered: $\gamma^{\delta}_1, \gamma^{\delta}_2, \ldots$ .

Each exploration process of course requires choosing a starting
and ending vertex and edge.
For domains of type~1, with a $\pm$ or $\mp$ boundary condition,
the choice is the natural one, explained before.

For a domain $D^{\delta}_k$ (used at the $k$th step) of type
other than~1, and therefore with a monochromatic boundary, the
starting and ending edges are chosen with a procedure that mimics
what is done in the continuum case. 
Once again, the exact procedure used to choose the pair of
points is not important, as long as they are not chosen too
close to each other.
This is clear in the discrete case because the procedure
that we are presenting is only ``discovering'' the cluster
boundaries.
In more precise terms, it is clear that one could couple
the processes obtained with different rules by means of the
same percolation configuration, thus obtaining exactly the
same cluster boundaries.

As in the continuum case, we can choose the following procedure.
(In Theorem~\ref{thm-convergence} we will slightly reorganize the
procedure by using a coupling to the continuum construction to
guarantee that the order of exploration of domains of the discrete
and continuum procedures match despite the rules for breaking ties.)
Given a domain $D$, $x$ and $y$ are chosen so that, of all pairs
$(u,v)$ of points in $\partial D$, they maximize $|\text{Re}(u-v)|$
if $|\text{Re}(u-v)| \geq |\text{Im}(u-v)|$, or else they maximize
$|\text{Im}(u-v)|$.
If the choice is not unique, to restrict the number of pairs
one looks at those pairs, among the ones already obtained, that
maximize the other of $\{ |\text{Re}(u-v)|, |\text{Im}(u-v)| \}$.
Notice that this leaves at most two pairs of points; if that's
the case, the pair that contains the point with minimal real
(and, if necessary, imaginary) part is chosen.

The procedure continues iteratively, with regions that
have monochromatic boundaries playing  the role played
in the first step by the unit disc.
Every time a region with $\pm$ boundary conditions is used,
a new loop, corresponding to the outer boundary contour of
a cluster, is formed by pasting together, as explained in
Section~\ref{single-loop}, the new exploration path and the
excursion containing the region where the last exploration
was carried out.
All the new regions created at a step when a loop is formed
have monochromatic boundary conditions.

\section{Main Technical Results} \label{main-results}

In this section we collect our main results about the Continuum
Nonsimple Loop process.
Before doing that, we state a precise version, called statement (S),
of convergence of exploration paths to $SLE_6$ that we will use in
the proofs of these results, presented in Section~\ref{proofs}.
Statement (S) is an immediate consequence of Theorem~5 of~\cite{cn2}.
The proof given in~\cite{cn2}, which relies among other things on the
result of Smirnov~\cite{smirnov} concerning convergence of
crossing probabilities to Cardy's formula~\cite{cardy,cardy2},
is an expanded and corrected version of Appendix~A of~\cite{cn1}.
We note that (S) is both more general and more special than the
convergence statements in~\cite{smirnov,smirnov-long} --- more
general in that the domain can vary with $\delta$ as $\delta \to 0$,
but more special in the restriction to Jordan domains.

Given a Jordan domain $D$ with two distinct points $a,b \in \partial D$
on its boundary, let $\mu_{D,a,b}$ denote the law of $\gamma_{D,a,b}$,
the trace of chordal $SLE_6$, and let $\mu^{\delta}_{D,a,b}$ denote
the law of the percolation exploration path $\gamma^{\delta}_{D,a,b}$.
Let $W$ be the space of continuous curves inside $D$ from $a$ to $b$.
We define $\rho(\mu_{D,a,b},\mu^{\delta}_{D,a,b}) \equiv
\inf \{\varepsilon>0 : \mu_{D,a,b}(U) \leq
\mu^{\delta}_{D,a,b}(\bigcup_{x \in U}B_{\text{d}}(x,\varepsilon))
+ \varepsilon \text{ for all Borel } U \subset W \}$
(where $B_{\text{d}}(x,\varepsilon)$ denotes the open ball of
radius $\varepsilon$ centered at $x$ in the metric~(\ref{distance}))
and denote by
$\text{d}_{\text{P}}(\mu_{D,a,b},\mu^{\delta}_{D,a,b}) \equiv
\max \{ \rho(\mu_{D,a,b},\mu^{\delta}_{D,a,b}), \rho(\mu^{\delta}_{D,a,b},\mu_{D,a,b}) \}$
the Prohorov distance; weak convergence is equivalent to
convergence in the Prohorov metric.
Statement (S) is the following; it is used in the proofs of all the results
of this section {\it except\/} for Lemmas~\ref{sub-conv}-\ref{boundaries}.
\begin{itemize}
\item[(S)] For Jordan domains, there is convergence in distribution
of the percolation exploration path to the trace of chordal $SLE_6$
that is \emph{locally uniform} in the shape of the boundary with respect
to the uniform metric on continuous curves~(\ref{distance}), and in the
location of the starting and ending points with respect to the Euclidean
metric; i.e., for $(D,a,b)$ a Jordan domain with distinct $a,b \in \partial D$,
$\forall \varepsilon>0$, $\exists \alpha_0=\alpha_0(\varepsilon)$ and
$\delta_0=\delta_0(\varepsilon)$ such that for all $(D',a',b')$ with
$D'$ Jordan and with
$\max{(\text{d}(\partial D, \partial D'),|a-a'|,|b-b'|) \leq \alpha_0}$
and $\delta \leq \delta_0$,
$\text{d}_{\text{P}}(\mu_{D',a',b'},\mu^{\delta}_{D',a',b'}) \leq \varepsilon$.
\end{itemize}

\subsection{Preliminary Results} \label{pre-res}

We first give some important results which
are needed in the proofs of the main theorems.
We start with two lemmas which are consequences
of~\cite{ab}, of standard bounds on the probability
of events corresponding to having a certain number of
monochromatic crossings of an annulus (see Lemma~5
of~\cite{ksz}, Appendix~A of~\cite{lsw5}, and also~\cite{ada}),
but which do \emph{not} depend on statement (S).

\begin{lemma} \label{sub-conv}
Let $\gamma^{\delta}_{{\mathbb D},-i,i}$ be the percolation
exploration path on the edges of $\delta {\cal H}$ inside
(the $\delta$-approximation of) $\mathbb D$ between
(the e-vertices closest to) $-i$ and $i$.
For any fixed point $z \in {\mathbb D}$, chosen independently
of $\gamma^{\delta}_{{\mathbb D},-i,i}$, as $\delta \to 0$,
$\gamma^{\delta}_{{\mathbb D},-i,i}$ and the boundary
$\partial {\mathbb D}^{\delta}_{-i,i}(z)$ of the domain
${\mathbb D}^{\delta}_{-i,i}(z)$ that contains $z$ jointly
have limits in distribution along subsequences of $\delta$
with respect to the uniform metric~(\ref{distance}) on
continuous curves.
Moreover, any subsequence limit of
$\partial {\mathbb D}^{\delta}_{-i,i}(z)$ is almost surely
a simple loop~\cite{ada}.
\end{lemma}


\begin{lemma} \label{boundaries}
Using the notation of Lemma~\ref{sub-conv},
let $\gamma_{{\mathbb D},-i,i}$ be the limit in distribution
of $\gamma^{\delta}_{{\mathbb D},-i,i}$ as $\delta \to 0$
along some convergent subsequence $\{ \delta_k \}$ and
$\partial {\mathbb D}_{-i,i}(z)$ the boundary of the domain
${\mathbb D}_{-i,i}(z)$ of ${\mathbb D} \setminus \gamma_{D,-i,i}[0,1]$
that contains $z$.
Then, as $k \to \infty$,
$(\gamma^{\delta_k}_{{\mathbb D},-i,i},\partial {\mathbb D}^{\delta_k}_{-i,i}(z))$
converges in distribution to
$(\gamma_{{\mathbb D},-i,i},\partial {\mathbb D}_{-i,i}(z))$.
\end{lemma}


The two lemmas above are important ingredients in the
proof of Theorem~\ref{thm-convergence} below.
The second one says that, for every subsequence limit, the discrete
boundaries converge to the boundaries of the domains generated by
the limiting continuous curve.
If we use statement (S), then the limit $\gamma_{{\mathbb D},-i,i}$
of $\gamma^{\delta_k}_{{\mathbb D},-i,i}$ is the trace of chordal
$SLE_6$ for every subsequence $\delta_k \downarrow 0$, and we can
use Lemmas~\ref{boundaries} and \ref{sub-conv} to deduce that all
the domains produced in the continuum construction are Jordan domains.
The key step in that direction is represented by the following
result, our proof of which relies on (S).

\begin{corollary} \label{jordan}
For any deterministic $z \in {\mathbb D}$,
the boundary $\partial {\mathbb D}_{-i,i}(z)$ of a domain
${\mathbb D}_{-i,i}(z)$ of the continuum construction is
almost surely a Jordan curve.
\end{corollary}

\noindent The corollary says that the domains that appear after
the first step of the continuum construction are Jordan domains.
The steps in the second stage of the continuum construction consist
of $SLE_6$ paths inside Jordan domains, and therefore Corollary~\ref{jordan},
combined with Riemann's mapping theorem and the conformal invariance
of $SLE_6$, implies that the domains produced during the second stage
are also Jordan.
By induction, we deduce that all the domains produced in the
continuum construction are Jordan domains.

We end this section with one more lemma which is another
key ingredient in the proof of Theorem~\ref{thm-convergence};
we remark that its proof requires (S) in a fundamental way.

\begin{lemma} \label{strong-smirnov}
Let $(D,a,b)$ denote a \emph{random} Jordan domain, with $a,b$ two
points on $\partial D$.
Let $\{ (D_k,a_k,b_k) \}_{k \in {\mathbb N}}, \, a_k,b_k \in \partial D_k$,
be a sequence of \emph{random} Jordan domains with points on their
boundaries such that, as $k \to \infty$, $(\partial D_k,a_k,b_k)$
converges in distribution to $(\partial D,a,b)$ with respect to the
uniform metric~(\ref{distance}) on continuous curves, and the Euclidean
metric on $(a,b)$.
For any sequence $\{ \delta_k \}_{k \in {\mathbb N}}$
with $\delta_k \downarrow 0$ as $k \to \infty$,
$\gamma^{\delta_k}_{D_k,a_k,b_k}$ converges in distribution to
$\gamma_{D,a,b}$ with respect to the uniform metric~(\ref{distance})
on continuous curves.
\end{lemma}

\subsection{Main Technical Theorems} \label{main-thms}

In this section we state the main technical theorems of this paper.
Our main results, presented in Section~\ref{quickresults}, are
consequences of these theorems.
The proofs of these theorems rely on statement (S).
As noted before, a detailed proof of statement (S) can be found
in~\cite{cn2}.

\begin{theorem} \label{thm-convergence} For
any $k \in {\mathbb N}$, the
first $k$ steps of (a suitably reorganized version of)
the full discrete construction inside the unit disc (of
Section~\ref{full}) converge, jointly in distribution,
to the first $k$ steps of the full continuum construction
inside the unit disc (of Section~\ref{unit-disc}).
Furthermore, the scaling limit of the full (original or
reorganized) discrete construction is the full continuum
construction.

Moreover, if for any fixed $\varepsilon>0$ we let
$K_{\delta}(\varepsilon)$ denote the number of steps
needed to find all the cluster boundaries of Euclidean
diameter larger than $\varepsilon$ in the discrete
construction, then $K_{\delta}(\varepsilon)$ is
bounded in probability as $\delta \to 0$; i.e.,
$\lim_{C \to \infty} \limsup_{\delta \to 0}
{\mathbb P}(K_{\delta}(\varepsilon) > C) = 0$.
This is so in both the original and reorganized
versions of the discrete construction.
\end{theorem}


The second part of Theorem~\ref{thm-convergence}
means that both versions of the discrete
construction used in the theorem find all large
contours in a number of steps which does not diverge
as $\delta \to 0$.
This, together with the first part of the same theorem,
implies that the continuum construction does indeed
describe all macroscopic contours contained inside
the unit disc (with blue boundary conditions) as
$\delta \to 0$.

The construction presented in Section~\ref{unit-disc}
can of course be repeated
for the disc ${\mathbb D}_R$ of radius $R$, for any $R$,
so we should take a ``thermodynamic limit'' by letting
$R \to \infty$.
In this way, we would eliminate the boundary (and
the boundary conditions) and obtain a process on
the whole plane.
Such an extension from the unit disc to the plane
is contained in the next theorem.

Let $P_R$ be the (limiting) distribution of the set of
curves (all continuum nonsimple loops) generated by the
continuum construction inside ${\mathbb D}_R$ (i.e., the
limiting measure, defined by the inductive construction,
on the complete separable metric space ${\Omega}_R$ of
collections of continuous curves in ${\mathbb D}_R$).

For a domain $D$, we denote by $I_D$ the mapping (on $\Omega$
or $\Omega_R$) in which all portions of curves that exit $D$
are removed.
When applied to a configuration of loops in the plane,
$I_D$ gives a set of curves which either start and end
at points on $\partial D$ or form closed loops completely
contained in $D$.
Let $\hat I_D$ be the same mapping lifted to the space
of probability measures on $\Omega$ or $\Omega_R$.

\begin{theorem} \label{thm-therm-lim}
There exists a unique probability measure $P$ on the space
$\Omega$ of collections of continuous curves in $\dot{\mathbb R}^2$
such that $P_R \to P$ as $R \to \infty$ in the sense that for every
bounded domain $D$, as $R \to \infty$, $\hat I_D P_R \to \hat I_D P$.
\end{theorem}

\begin{remark}
We we will generally take monochromatic blue boundary
conditions on the disc ${\mathbb D}_R$ of radius $R$,
but this arbitrary choice does not affect the results.
\end{remark}


The next theorem states a conformal invariance property
of the Continuum Nonsimple Loop processes of Theorem~\ref{main-result1}.


\begin{theorem} \label{thm-conf-inv}
Given two disjoint discs, $D_1$ and $D_2$, let $\lambda_1$
(respectively, $\lambda_2$) be the smallest loop from the
Continuum Nonsimple Loop process X that surrounds $D_1$
(resp., $D_2$) and let $\tilde D_1$ (resp., $\tilde D_2$) be
the connected component of ${\mathbb R}^2 \setminus \lambda_1$
(resp., ${\mathbb R}^2 \setminus \lambda_2$) that contains
$D_1$ (resp., $D_2$).
Assume that $\tilde D_1$ and $\tilde D_2$ are disjoint and
let $P_{\tilde D_i}$, $i=1,2$, denote the distribution of
the loops
inside $\tilde D_i$.
Then, conditioned on $\tilde D_1$ and $\tilde D_2$, the
configurations inside $\tilde D_1$ and $\tilde D_2$ are
independent and moreover $P_{\tilde D_2} = f * P_{\tilde D_1}$
(here $f * P_{\tilde D_1}$ denotes the probability distribution
of the loop process $f(X')$ when $X'$ is distributed by $P_{\tilde D_1}$),
where $f : \tilde D_1 \to \tilde D_2$ is a conformal homeomorphism
from $\tilde D_1$ onto $\tilde D_2$.
\end{theorem}

We remark that the result is still valid (without the independence)
even if $\tilde D_1$ and $\tilde D_2$ are not disjoint,
but for simplicity we do not consider that case.

To conclude this section, we show how to recover chordal $SLE_6$
from the Continuum Nonsimple Loop process, i.e., given a (deterministic)
Jordan domain $D$ with two boundary points $a$ and $b$, we give a
construction that uses the continuum nonsimple loops of $P$ to generate
a process distributed like chordal $SLE_6$ inside $D$ from $a$ to $b$.

Remember, first of all, that each continuum nonsimple loop has
either a clockwise or counterclockwise direction, with the set
of all loops surrounding any deterministic point alternating
in direction.
For convenience, let us suppose that $a$ is at the ``bottom"
and $b$ is at the ``top" of $D$ so that the boundary is divided
into a left and right part by these two points.
Fix $\varepsilon>0$ and call $LR(\varepsilon)$ the set of all
the directed segments of loops that connect from the left to
the right part of the boundary touching $\partial D$ at a distance
larger than $\varepsilon$ from both $a$ and $b$, and $RL(\varepsilon)$
the analogous set of directed segments from the right to the left
portion of $\partial D$.
For a fixed $\varepsilon>0$, there is only a finite number
of such segments, and, if they are ordered moving along the
left boundary of $D$ from $a$ to $b$, they alternate in direction
(i.e., a segment in $LR(\varepsilon)$ is followed by one in
$RL(\varepsilon)$ and so on).

Between a segment in $RL(\varepsilon)$ and the next segment
in $LR(\varepsilon)$, there are countably many portions of
loops intersecting $D$ which start and end on $\partial D$
and are maximal in the sense that they are not contained
inside any other portion of loop of the same type; they all
have counterclockwise direction and can be used to make a
``bridge'' between the right-to-left segment and the next
one (in $LR(\varepsilon)$).
This is done by pasting the portions of loops together with
the help of points in $\partial D$ and a limit procedure
to produce a connected (nonsimple) path.

If we do this for each pair of successive segments on both
sides of the boundary of $D$, we get a path that connects
two points on $\partial D$.
By letting $\varepsilon \to 0$ and taking the limit of this
procedure, since almost surely $a$ and $b$ are surrounded
by an infinite family of nested loops with diameters going
to zero, we obtain a path that connects $a$ with $b$;
this path is distributed as chordal $SLE_6$ inside $D$
from $a$ to $b$.
The last claim follows from considering the analogous
procedure for percolation on the discrete lattice
$\delta {\cal H}$, using segments of boundaries.
It is easy to see that in the discrete case this
procedure produces exactly the same path as the
percolation exploration process.
By Theorems \ref{main-result1} and~\ref{main-result2},
the scaling limit of this discrete procedure is the
continuum one described above, therefore the claim
follows from (S).

%
%

\section{Proofs} \label{proofs}

In this section we present the proofs of the results
stated in Sections~\ref{quickresults} and~\ref{main-results}.
In order to do that, we will use the following lemma.
\begin{lemma} \label{arms}
Let ${\cal A}^{\delta}(v;\varepsilon,\varepsilon')$ be the event that
the annulus $B(v,\varepsilon) \setminus B(v,\varepsilon')$ centered at
$v \in {\mathbb D}$ contains six disjoint monochromatic crossings, not
all of the same color, and let ${\cal B}^{\delta}(v;\varepsilon,\varepsilon')$
be the event, for some $v \in \partial{\mathbb D}$, that
${\mathbb D} \cap \{B(v,\varepsilon) \setminus B(v,\varepsilon')\}$
contains three disjoint monochromatic crossings,
not all of the same color.
Then, for any $\varepsilon>0$, 
\begin{equation} \label{six}
\lim_{\varepsilon' \to 0} \, \limsup_{\delta \to 0} \,
{\mathbb P}(\bigcup_{v \in {\mathbb D}} {\cal A}^{\delta}(v;\varepsilon,\varepsilon')) = 0
\end{equation}
and
\begin{equation} \label{three}
\lim_{\varepsilon' \to 0} \, \limsup_{\delta \to 0} \,
{\mathbb P}(\bigcup_{v \in \partial{\mathbb D}} {\cal B}^{\delta}(v;\varepsilon,\varepsilon')) = 0.
\end{equation}
\end{lemma}

\noindent {\bf Proof.} We know from~\cite{ksz} that
there exist $c_1 < \infty$ and $\alpha>0$ so that for
$\varepsilon_2<\varepsilon_1$, and $\delta$ small enough (in particular, $\delta<\varepsilon_2$),
\begin{equation} \label{ksz}
{\mathbb P}({\cal A}^{\delta}(v;\varepsilon_1,\varepsilon_2))
\leq c_1 \, \left(\frac{\varepsilon_2}{\varepsilon_1}\right)^{2+\alpha}
\end{equation}
for any $v \in {\mathbb R}^2$.
If we cover $\mathbb D$ with $N_{\varepsilon'}$ balls of radius
$\varepsilon'$ centered at $\{ v_j \}_{j \in {\cal N}_{\varepsilon'}}$,
we have that, for $\varepsilon' < \varepsilon /6$ and $\delta$ small enough,
\begin{equation}
{\mathbb P}(\bigcup_{v \in \mathbb D} {\cal A}^{\delta}(v;\varepsilon,\varepsilon')) \leq
{\mathbb P}(\bigcup_{j \in {\cal N}_{\varepsilon'}} {\cal A}^{\delta}(v_j;\varepsilon/2,3 \, \varepsilon'))
\leq 6^{2+\alpha} \, c_1 \, N_{\varepsilon'} \left(\frac{\varepsilon'}{\varepsilon}\right)^{2+\alpha},
\end{equation}
where the first inequality follows from the observations that for any
$v \in \mathbb D$, $B(v,\varepsilon') \subset B(v_j,3 \, \varepsilon')$
and $B(v_j,\varepsilon/2) \subset B(v_j,\varepsilon - \varepsilon') \subset B(v,\varepsilon)$
for some $j \in {\cal N}_{\varepsilon'}$, and the second inequality uses~(\ref{ksz}).
Using the fact that $N_{\varepsilon'}$ is $O(\frac{1}{\varepsilon'})^2$,
we can let first $\delta \to 0$ 
and then
$\varepsilon' \to 0$ to obtain~(\ref{six}).

We also know, as a consequence of Lemma~5 of~\cite{ksz} or as
proved in Appendix A of~\cite{lsw5},
that for any $v \in {\mathbb R}$, the probability that the semi-annulus
${\mathbb H} \cap \{B(v,\varepsilon_1) \setminus B(v,\varepsilon_2)\}$
contains three disjoint monochromatic crossings,
not all of the same color, is bounded above by
$c_2 \, (\varepsilon_2/\varepsilon_1)^{1+\beta}$ for some $c_2<\infty$
and $\beta>0$.
(We remark that the result still applies when ${\mathbb H}$
is replaced by any other half-plane.)
Since the unit disc is a convex subset of the half-plane $\{ x+iy:y>-1 \}$
and therefore the intersection of an annulus centered at $-i$ with
the unit disc $\mathbb D$ is a subset of the intersection of the same
annulus with the half-plane $\{ x+iy:y>-1 \}$, we can use that bound to
conclude that for $v=-i$, and in fact
for any $v \in \partial{\mathbb D}$,
there exists a constant $c_2 < \infty$ such that
\begin{equation} \label{lsw}
{\mathbb P}({\cal B}^{\delta}(v;\varepsilon_1,\varepsilon_2))
\leq c_2 \, \left(\frac{\varepsilon_2}{\varepsilon_1}\right)^{1+\beta}
\end{equation}
for some $\beta>0$.
We can then use similar arguments to those above, together with~(\ref{lsw}),
to obtain~(\ref{three}) and conclude the proof.~\fbox{} \\

\noindent {\bf Proof of Lemma~\ref{sub-conv}.}
The first part of the lemma is a direct consequence
of~\cite{ab}; it is enough to notice that the (random)
polygonal curves $\gamma^{\delta}_{{\mathbb D},-i,i}$
and $\partial {\mathbb D}^{\delta}_{-i,i}(z)$
satisfy the conditions in~\cite{ab} and thus have a
scaling limit in terms of continuous curves, at least
along subsequences of $\delta$.

To prove the second part, we use a standard percolation
bound (see Lemma~5 of~\cite{ksz}) 
to show that, in the limit $\delta \to 0$, the loop
$\partial {\mathbb D}^{\delta}_{-i,i}(z)$ does not collapse
on itself but remains a  simple loop.

Let us assume that this is not the case and that
the limit $\tilde\gamma$ of
$\partial {\mathbb D}^{\delta_k}_{-i,i}(z)$ along some
subsequence $\{ \delta_k \}_{k \in {\mathbb N}}$ touches
itself, i.e., $\tilde\gamma(t_0)=\tilde\gamma(t_1)$ for
$t_0 \neq t_1$ with positive probability.
If that happens, we can take $\varepsilon>\varepsilon'>0$
small enough so that the annulus
$B(\tilde\gamma(t_1),\varepsilon) \setminus B(\tilde\gamma(t_1),\varepsilon')$
is crossed at least four times by $\tilde\gamma$
(here $B(u,r)$ is the ball of radius $r$ centered at $u$).

Because of the choice of topology, the convergence in
distribution of $\partial {\mathbb D}^{\delta_k}_{-i,i}(z)$
to $\tilde\gamma$ implies that we can find coupled versions of
$\partial {\mathbb D}^{\delta_k}_{-i,i}(z)$ and $\tilde\gamma$
on some probability space $(\Omega',{\cal B}',{\mathbb P}')$ such that
$\text{d}(\partial {\mathbb D}^{\delta}_{-i,i}(z),\tilde\gamma) \to 0$,
for all $\omega' \in \Omega'$ as $k \to \infty$ (see, for example,
Corollary~1 of~\cite{billingsley1}).

\begin{figure}[!ht]
\begin{center}
\includegraphics[width=8cm]{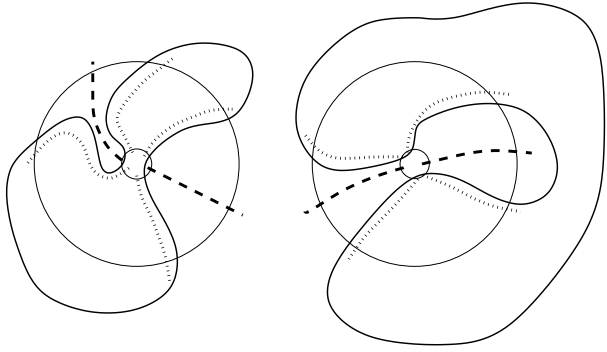}
\caption{Schematic diagrams representing four blue (dotted
in the figure) and two yellow (dashed
in the figure) crossings of an annulus produced
by having four crossings of the same annulus by a boundary
(the solid loops).}
\label{fig-lemma5-1}
\end{center}
\end{figure}

Using this coupling, we can choose $k$ large enough (depending
on $\omega'$) so that $\partial {\mathbb D}^{\delta_k}_{-i,i}(z)$
stays in an $\varepsilon'/2$-neighborhood
${\cal N}(\tilde\gamma,\varepsilon'/2) \equiv \bigcup_{u \in \tilde\gamma} B(u,\varepsilon'/2)$
of $\tilde\gamma$.
This event however would correspond to (at least) four
paths of one color (corresponding to the four crossings
by $\Delta {\mathbb D}^{\delta_k}_{-i,i}(z)$, which
shadows
$\partial {\mathbb D}^{\delta_k}_{-i,i}(z)$) and
two of the other color (belonging to percolation clusters
adjacent to the cluster of $\Delta {\mathbb D}^{\delta_k}_{-i,i}(z)$,
and of the opposite color), of the annulus
$B(\tilde\gamma(t_1),\varepsilon-\varepsilon'/2)
\setminus B(\tilde\gamma(t_1),3 \, \varepsilon'/2)$ (see,
for example,~\cite{ada} --- see also Figure~\ref{fig-lemma5-1}).
As $\delta_k \to 0$, we can let $\varepsilon' \to 0$,
in which case the probability of seeing the event
just described somewhere inside $\mathbb D$ goes to
zero by an application of Lemma~\ref{arms}, leading
to a contradiction. \fbox{} \\

In order to prove Lemma~\ref{boundaries}, we will use the following result.
\begin{lemma} \label{claim-c}
For two (deterministic) points $u,v \in {\mathbb D}$,
the probability that ${\mathbb D}_{-i,i}(u) = {\mathbb D}_{-i,i}(v)$ but
${\mathbb D}^{\delta}_{-i,i}(u) \neq {\mathbb D}^{\delta}_{-i,i}(v)$
or vice versa goes to zero as $\delta \to 0$.
\end{lemma}

\noindent {\bf Proof.}
Let $\{ \delta_k \}_{k \in {\mathbb N}}$ be a convergent
subsequence for $\gamma^{\delta}_{{\mathbb D},-i,i}$
and let $\gamma \equiv \gamma_{{\mathbb D},-i,i}$ be the
limit in distribution of $\gamma^{\delta_k}_{{\mathbb D},-i,i}$
as $k \to \infty$.
For simplicity of notation, in the rest of the proof we will drop
the $k$ and write $\delta$ instead of $\delta_k$.
Because of the choice of topology, the convergence in distribution
of $\gamma^{\delta} \equiv \gamma^{\delta}_{{\mathbb D},-i,i}$ to
$\gamma$ implies that we can find coupled versions of $\gamma^{\delta}$
and $\gamma$ on some probability space $(\Omega',{\cal B}',{\mathbb P}')$
such that $\text{d}(\gamma^{\delta}(\omega'),\gamma(\omega')) \to 0$,
for all $\omega'$ as $k \to \infty$ (see, for example, Corollary~1
of~\cite{billingsley1}).

Using this coupling, we first consider the case of $u,v$
such that ${\mathbb D}_{-i,i}(u) = {\mathbb D}_{-i,i}(v)$ but
${\mathbb D}^{\delta}_{-i,i}(u) \neq {\mathbb D}^{\delta}_{-i,i}(v)$.
Since ${\mathbb D}_{-i,i}(u)$ is an open subset of $\mathbb C$,
there exists a continuous curve $\gamma_{u,v}$ joining
$u$ and $v$ and a constant $\varepsilon>0$ such that
the $\varepsilon$-neighborhood
${\cal N}(\gamma_{u,v},\varepsilon)$ of the curve is contained
in ${\mathbb D}_{-i,i}(u)$, which implies that $\gamma$ does
not intersect ${\cal N}(\gamma_{u,v},\varepsilon)$.
Now, if $\gamma^{\delta}$ does not intersect
${\cal N}(\gamma_{u,v},\varepsilon/2)$, for
$\delta$ small enough, then there is a $\cal T$-path
$\pi$ of unexplored hexagons connecting the hexagon
that contains $u$ with the hexagon that contains $v$,
and we conclude that
${\mathbb D}^{\delta}_{-i,i}(u) = {\mathbb D}^{\delta}_{-i,i}(v)$.

This shows that the event that ${\mathbb D}_{-i,i}(u) = {\mathbb D}_{-i,i}(v)$
but ${\mathbb D}^{\delta}_{-i,i}(u) \neq {\mathbb D}^{\delta}_{-i,i}(v)$
implies the existence of a curve $\gamma_{u,v}$ whose
$\varepsilon$-neighborhood ${\cal N}(\gamma_{u,v},\varepsilon)$
is not intersected by $\gamma$ but whose $\varepsilon/2$-neighborhood
${\cal N}(\gamma_{u,v},\varepsilon/2)$ is intersected by $\gamma^{\delta}$.
This implies that $\forall u,v \in {\mathbb D}$, $\exists \varepsilon>0$ such
that ${\mathbb P}'({\mathbb D}_{-i,i}(u) = {\mathbb D}_{-i,i}(v) \text{ but }
{\mathbb D}^{\delta}_{-i,i}(u) \neq {\mathbb D}^{\delta}_{-i,i}(v)) \leq
{\mathbb P}'(\text{d}(\gamma^{\delta},\gamma) \geq \varepsilon/2)$.
But the right hand side goes to zero for every $\varepsilon>0$ as
$\delta \to 0$, which concludes the proof of one direction of the claim.

To prove the other direction, we consider two points
$u,v \in D$ such that $D_{-i,i}(u) \neq D_{-i,i}(v)$
but $D^{\delta}_{-i,i}(u) = D^{\delta}_{-i,i}(v)$.
Assume that $u$ is trapped before $v$ by $\gamma$
and suppose for the moment that ${\mathbb D}_{-i,i}(u)$
is a domain of type 3 or 4; the case of a domain of
type 1 or 2 is analogous and will be treated later.
Let $t_1$ be the first time $u$ is trapped by
$\gamma$ with $\gamma(t_0)=\gamma(t_1)$ the double
point of $\gamma$ where the domain ${\mathbb D}_{-i,i}(u)$
containing $u$ is ``sealed off.''
At time $t_1$, a new domain containing $u$ is
created and $v$ is disconnected from $u$.

Choose $\varepsilon>0$ small enough so that neither $u$ nor
$v$ is contained in the ball $B(\gamma(t_1),\varepsilon)$
of radius $\varepsilon$ centered at $\gamma(t_1)$, nor in
the $\varepsilon$-neighborhood
${\cal N}(\gamma[t_0,t_1],\varepsilon)$ of the portion
of $\gamma$ which surrounds $u$. 
Then it follows from the coupling that, for $\delta$
small enough, there are appropriate parameterizations
of $\gamma$ and $\gamma^{\delta}$ such that the portion
$\gamma^{\delta}[t_0,t_1]$ of $\gamma^{\delta}(t)$ is
inside ${\cal N}(\gamma[t_0,t_1],\varepsilon)$, and
$\gamma^{\delta}(t_0)$ and $\gamma^{\delta}(t_1)$ are
contained in $B(\gamma(t_1),\varepsilon)$.

\begin{figure}[!ht]
\begin{center}
\includegraphics[width=8cm]{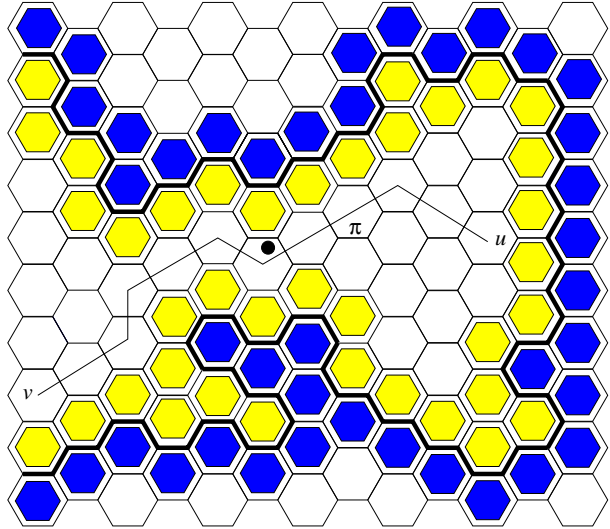}
\caption{Example of a $\cal T$-path $\pi$ of unexplored
hexagons from $u$ to $v$ having to go through a ``bottleneck"
due to the fact that the exploration path (heavy line) comes
close to itself.
An approximate location of the continuum double point at
$\gamma(t_0) = \gamma(t_1)$ is indicated by the small disc
in one of the hexagons in the bottleneck area.}
\label{fig-lemma6-2}
\end{center}
\end{figure}

For $u$ and $v$ to be contained in the same domain
in the discrete construction, there must be a
$\cal T$-path $\pi$ of unexplored hexagons
connecting the hexagon that contains $u$ to the
hexagon that contains $v$.
>From what we said in the previous paragraph, any
such $\cal T$-path connecting $u$ and $v$
would have to go though a ``bottleneck'' in
$B(\gamma(t_1),\varepsilon)$ (see Figure~\ref{fig-lemma6-2}).

Assume now, for concreteness but without loss of
generality, that ${\mathbb D}_{-i,i}(u)$ is a domain
of type~3, which means that $\gamma$ winds around
$u$ counterclockwise, and consider the hexagons
to the ``left" of $\gamma^{\delta}[t_0,t_1]$
(these are all lightly shaded in Figure~\ref{fig-lemma6-2}).
Those hexagons form a ``quasi-loop'' around $u$
since they wind around it (counterclockwise) and
the first and last hexagons are both contained in
$B(\gamma(t_1),\varepsilon)$.
The hexagons to the left of $\gamma^{\delta}[t_0,t_1]$
belong to the set $\Gamma_Y(\gamma^{\delta})$, which
can be seen as a (nonsimple) path by connecting the
centers of the hexagons in $\Gamma_Y(\gamma^{\delta})$
by straight segments.
Such a path shadows $\gamma^{\delta}$, with the difference
that it can have double (or even triple) points, since
the same hexagon can be visited more than once.
Consider $\Gamma_Y(\gamma^{\delta})$ as a path
$\hat\gamma^{\delta}$ with a given parametrization
$\hat\gamma^{\delta}(t)$, chosen so that
$\hat\gamma^{\delta}(t)$ is inside
$B(\gamma(t_1),\varepsilon)$ when $\gamma^{\delta}(t)$ is,
and it winds around $u$ together with $\gamma^{\delta}(t)$.

Now suppose that there were two times, $\hat t_0$ and
$\hat t_1$, such that $\hat\gamma^{\delta}(\hat t_1)
= \hat\gamma^{\delta}(\hat t_0) \in B(\gamma(t_1),\varepsilon)$
and $\hat\gamma^{\delta}[\hat t_0,\hat t_1]$ winds
around $u$.
This would imply that the ``quasi-loop'' of explored
yellow hexagons around $u$ is actually completed, and
that $D^{\delta}_{a,b}(v) \neq D^{\delta}_{a,b}(u)$.
Thus, for $u$ and $v$ to belong to the same discrete
domain, this cannot happen.

For any $0<\varepsilon'<\varepsilon$, if we take $\delta$
small enough, $\hat\gamma^{\delta}$ will be contained
inside ${\cal N}(\gamma,\varepsilon')$, due to the coupling.
Following the considerations above, the fact that $u$
and $v$ belong to the same domain in the discrete
construction but to different domains in the continuum
construction implies, for $\delta$ small enough, that
there are four disjoint yellow $\cal T$-paths
crossing the annulus $B(\gamma(t_1),\varepsilon)
\setminus B(\gamma(t_1),\varepsilon')$ (the paths
have to be disjoint because, as we said,
$\hat\gamma^{\delta}$ cannot, when coming back to
$B(\gamma(t_1),\varepsilon)$ after winding around
$u$, touch itself inside $B(\gamma(t_1),\varepsilon)$).
Since $B(\gamma(t_1),\varepsilon) \setminus
B(\gamma(t_1),\varepsilon')$ is also crossed
by at least two blue $\cal T$-paths from
$\Gamma_B(\gamma^{\delta})$, there is a total
of at least six $\cal T$-paths, not all of
the same color, crossing the annulus
$B(\gamma(t_1),\varepsilon) \setminus B(\gamma(t_1),\varepsilon')$.
We can then use Lemma~\ref{arms} to conclude that,
if we keep $\varepsilon$ fixed and let $\delta \to 0$
and $\varepsilon' \to 0$, the probability to see such
an event anywhere in $\mathbb D$ goes to zero.



In the case in which $u$ belongs to a domain of type~1
or 2, let $\cal E$ be the excursion that traps $u$
and $\gamma(t_0) \in \partial {\mathbb D}$ be the point
on the boundary of $\mathbb D$ where $\cal E$ starts and
$\gamma(t_1) \in \partial {\mathbb D}$ the point where
it ends.
Choose $\varepsilon>0$ small enough so that neither $u$
nor $v$ is contained in the balls
$B(\gamma(t_0),\varepsilon)$ and $B(\gamma(t_1),\varepsilon)$
of radius $\varepsilon$ centered at $\gamma(t_0)$ and
$\gamma(t_1)$, nor in the $\varepsilon$-neighborhood
${\cal N}({\cal E},\varepsilon)$ of the excursion $\cal E$.
Because of the coupling, for $\delta$ small enough
(depending on $\varepsilon$), $\gamma^{\delta}$ shadows
$\gamma$ along $\cal E$, staying within
${\cal N}({\cal E},\varepsilon)$.
If this is the case, any $\cal T$-path of unexplored
hexagons connecting the hexagon that contains $u$ with
the hexagon that contains $v$ would have to go through
one of two ``bottlenecks,'' one contained in
$B(\gamma(t_0),\varepsilon)$ and the other in
$B(\gamma(t_1),\varepsilon)$.

Assume for concreteness (but without loss
of generality) that $u$ is in a domain of type 1,
which means that $\gamma$ winds around $u$
counterclockwise.
If we parameterize $\gamma$ and $\gamma^{\delta}$ so that
$\gamma^{\delta}(t_0) \in B(\gamma(t_0),\varepsilon)$
and $\gamma^{\delta}(t_1) \in B(\gamma(t_1),\varepsilon)$,
$\gamma^{\delta}[t_0,t_1]$ forms a ``quasi-excursion''
around $u$ since it winds around it (counterclockwise)
and it starts inside $B_{\varepsilon}(\gamma(t_0))$ and
ends inside $B_{\varepsilon}(\gamma(t_1))$.
Notice that if $\gamma^{\delta}$ touched
$\partial {\mathbb D}^{\delta}$, inside both
$B_{\varepsilon}(\gamma(t_0))$ and
$B_{\varepsilon}(\gamma(t_1))$, this would imply
that the ``quasi-excursion'' is a real excursion and
that $D^{\delta}_{a,b}(v) \neq D^{\delta}_{a,b}(u)$.

For any $0<\varepsilon'<\varepsilon$, if we take $\delta$
small enough, $\gamma^{\delta}$ will be contained inside
${\cal N}(\gamma,\varepsilon')$, due to the coupling.
Therefore, the fact that
${\mathbb D}^{\delta}_{a,b}(v) = {\mathbb D}^{\delta}_{a,b}(u)$
implies, with probability going to one as
$\delta \to 0$, that for $\varepsilon>0$ fixed
and any $0<\varepsilon'<\varepsilon$, $\gamma^{\delta}$
enters the ball $B(\gamma(t_i),\varepsilon')$
and does not touch $\partial {\mathbb D}^{\delta}$
inside the larger ball $B(\gamma(t_i),\varepsilon)$,
for $i=0$ or $1$.
This is equivalent to having at least two yellow and one
blue $\cal T$-paths (contained in ${\mathbb D}^{\delta}$)
crossing the annulus
$B(\gamma(t_i),\varepsilon) \setminus B(\gamma(t_i),\varepsilon')$.
As $\delta \to 0$, we can let $\varepsilon'$ go to zero
(keeping $\varepsilon$ fixed) and use Lemma~\ref{arms}
to conclude that the probability that such an event occurs
anywhere on the boundary of the unit disc goes to zero.

We have shown that, for two fixed points $u,v \in {\mathbb D}$,
having ${\mathbb D}_{-i,i}(u) \neq {\mathbb D}_{-i,i}(v)$ but
${\mathbb D}^{\delta}_{-i,i}(u) = {\mathbb D}^{\delta}_{-i,i}(v)$
or vice versa implies the occurrence of an event whose
probability goes to zero as $\delta \to 0$, and the
proof of the lemma is concluded.~\fbox{} \\

\noindent {\bf Proof of Lemma~\ref{boundaries}.}
As in the proof of Lemma~\ref{claim-c}, we let
$\{ \delta_k \}_{k \in {\mathbb N}}$ be a convergent
subsequence for $\gamma^{\delta}_{{\mathbb D},-i,i}$
and let $\gamma \equiv \gamma_{{\mathbb D},-i,i}$ be
the limit in distribution of $\gamma^{\delta_k}_{{\mathbb D},-i,i}$
as $k \to \infty$, and in the rest of the proof consider coupled
versions of $\gamma^{\delta_k} \equiv \gamma^{\delta_k}_{{\mathbb D},-i,i}$
and $\gamma$.

Let us introduce the Hausdorff distance $\text{d}_{\text{H}}(A,B)$
between two closed nonempty subsets of $\overline{\mathbb D}$:
\begin{equation} \label{hausdorff-dist}
\text{d}_{\text{H}}(A,B) \equiv \inf \{ \ell \geq 0 : B \subset
\cup_{a \in A} B(a,\ell), \, A \subset
\cup_{b \in B} B(b,\ell) \}.
\end{equation}
With this metric, the collection of closed subsets
of $\overline{\mathbb D}$ is a compact space.
We will next prove that
$\partial {\mathbb D}^{\delta_k}_{-i,i}(z)$ converges
in distribution to $\partial {\mathbb D}_{-i,i}(z)$
as $\delta_k \to 0$, in the topology induced
by~(\ref{hausdorff-dist}).
(Notice that the coupling between $\gamma^{\delta_k}$
and $\gamma$ provides a coupling between
$\partial {\mathbb D}^{\delta_k}_{-i,i}(z)$ and
$\partial {\mathbb D}_{-i,i}(z)$, seen as
boundaries of domains produced by the two paths.)


We will now use Lemma~\ref{sub-conv} and take a further subsequence
$k_n$ of the $\delta_k$'s that for simplicity of notation we denote
by $\{ \delta_n \}_{n \in {\mathbb N}}$ such that, as $n \to \infty$,
$\{ \gamma^{\delta_n},\partial {\mathbb D}^{\delta_n}_{-i,i}(z) \}$
converge jointly in distribution to $\{ \gamma,\tilde\gamma \}$,
where $\tilde\gamma$ is a simple loop.
For any $\varepsilon>0$, since $\tilde\gamma$ is a compact set,
we can find a covering of $\tilde \gamma$ by a finite number
of balls of radius $\varepsilon/2$ centered at points on $\tilde\gamma$.
Each ball contains both points in the interior
$\text{int}(\tilde\gamma)$ of $\tilde\gamma$ and in the exterior
$\text{ext}(\tilde\gamma)$ of $\tilde \gamma$, and we can choose
(independently of $n$) one point from $\text{int}(\tilde\gamma)$
and one from $\text{ext}(\tilde\gamma)$ inside each ball.

Once again, the convergence in distribution of
$\partial {\mathbb D}^{\delta_n}_{-i,i}(z)$
to $\tilde\gamma$ implies the existence of a coupling
such that, for $n$ large enough, the selected points
that are in $\text{int}(\tilde\gamma)$ are contained
in ${\mathbb D}^{\delta_n}_{-i,i}(z)$, and those that
are in $\text{ext}(\tilde\gamma)$ are contained in the
complement of $\overline{{\mathbb D}^{\delta_n}_{-i,i}(z)}$.
But by Lemma~\ref{claim-c}, each one of the selected points
that is contained in ${\mathbb D}^{\delta_n}_{-i,i}(z)$ is
also contained in ${\mathbb D}_{-i,i}(z)$ with probability
going to $1$ as $n \to \infty$; analogously, each one of the
selected points contained in the complement of
$\overline{{\mathbb D}^{\delta_n}_{-i,i}(z)}$ is also contained
in the complement of $\overline{{\mathbb D}_{-i,i}(z)}$ with
probability going to $1$ as $n \to \infty$.
This implies that $\partial {\mathbb D}_{-i,i}(z)$ crosses each
one of the balls in the covering of $\tilde\gamma$, and therefore
$\tilde\gamma \subset \cup_ {u \in \partial {\mathbb D}_{-i,i}(z)} B(u,\varepsilon)$.
>From this and the coupling between
$\partial {\mathbb D}^{\delta_n}_{-i,i}(z)$ and $\tilde\gamma$,
it follows immediately that, for $n$ large enough,
$\partial {\mathbb D}^{\delta_n}_{-i,i}(z) \subset
\cup_ {u \in \partial {\mathbb D}_{-i,i}(z)} B(u,\varepsilon)$
with probability close to one.

A similar argument (analogous to the previous one but simpler,
since it does not require the use of $\tilde\gamma$), with the roles
of ${\mathbb D}^{\delta_n}_{-i,i}(z)$ and ${\mathbb D}_{-i,i}(z)$
inverted, shows that $\partial {\mathbb D}_{-i,i}(z) \subset
\cup_ {u \in \partial {\mathbb D}^{\delta_n}_{-i,i}(z)} B(u,\varepsilon)$
with probability going to $1$ as $n \to \infty$.
Therefore, for all $\varepsilon>0$,
${\mathbb P}(\text{d}_{\text{H}}(\partial {\mathbb D}^{\delta_n}_{-i,i}(z),
\partial {\mathbb D}_{-i,i}(z)) > \varepsilon) \to 0$ as $n \to \infty$,
which implies convergence in distribution of
$\partial {\mathbb D}^{\delta_n}_{-i,i}(z)$ to
$\partial {\mathbb D}_{-i,i}(z)$, as $\delta_n \to 0$,
in the topology induced by~(\ref{hausdorff-dist}).
But Lemma~\ref{sub-conv} implies that
$\partial {\mathbb D}^{\delta_n}_{-i,i}(z)$
converges in distribution (using~(\ref{distance})) to a
simple loop, therefore $\partial {\mathbb D}_{-i,i}(z)$
must also be a simple loop; and we have convergence in
the topology induced by~(\ref{distance}).

It is also clear that the argument above is independent
of the subsequence $\{ \delta_n \}$
(and of the original subsequence $\{ \delta_k \}$), so
the limit of $\partial {\mathbb D}^{\delta}_{-i,i}(z)$ is
unique and coincides with $\partial {\mathbb D}_{-i,i}(z)$.
Hence, we have convergence in distribution of
$\partial {\mathbb D}^{\delta}_{-i,i}(z)$ to
$\partial {\mathbb D}_{-i,i}(z)$, as $\delta \to 0$,
in the topology induced by~(\ref{distance}), and
indeed joint convergence of
$(\gamma^{\delta},\partial {\mathbb D}^{\delta}_{-i,i}(z))$
to $(\gamma,\partial {\mathbb D}_{-i,i}(z))$. \fbox{} \\

\noindent {\bf Proof of Corollary~\ref{jordan}.}
The corollary follows immediately from Lemma~\ref{sub-conv}
and Lemma~\ref{boundaries}, as already seen in the proof
of Lemma~\ref{boundaries}. \fbox {} \\

\noindent {\bf Proof of Lemma~\ref{strong-smirnov}.}
First of all recall that the convergence of $(\partial D_k,a_k,b_k)$
to $(\partial D,a,b)$ in distribution implies the existence of
coupled versions of $(\partial D_k,a_k,b_k)$ and $(\partial D,a,b)$
on some probability space $(\Omega',{\cal B}',{\mathbb P}')$ such that
$\text{d}(\partial D(\omega'),\partial D_k(\omega')) \to 0$,
$a_k(\omega') \to a(\omega')$, $b_k(\omega') \to b(\omega')$
for all $\omega'$ as $k \to \infty$ (see, for example,
Corollary~1 of~\cite{billingsley1}).
This immediately implies that the conditions to apply Rad\'o's theorem
(see Theorem~\ref{rado-thm} of Section~\ref{rado}) are satisfied.
%
%
Let $f_k$ be the conformal map that takes the unit disc $\mathbb D$
onto $D_k$ with $f_k(0)=0$ and $f'_k(0)>0$, and let $f$ be the
conformal map from $\mathbb D$ onto $D$ with $f(0)=0$ and
$f'(0)>0$.
Then, by Theorem~\ref{rado-thm}, $f_k$ converges to $f$
uniformly in $\overline{\mathbb D}$, as $k \to \infty$.

Let $\gamma$ (resp., $\gamma_k$) be the chordal $SLE_6$
inside $D$ (resp., $D_k$) from $a$ to $b$ (resp., from $a_k$
to $b_k$), $\tilde \gamma = f^{-1}(\gamma)$, $\tilde a = f^{-1}(a)$,
$\tilde b = f^{-1}(b)$, and $\tilde \gamma_k = f_k^{-1}(\gamma_k)$,
$\tilde a_k = f_k^{-1}(a_k)$, $\tilde b_k = f_k^{-1}(b_k)$.
We note that, because of the conformal invariance of chordal
$SLE_6$, $\tilde\gamma$ (resp., $\tilde\gamma_k$) is distributed
as chordal $SLE_6$ in $\mathbb D$ from $\tilde a$ to $\tilde b$
(resp., from $\tilde a_k$ to $\tilde b_k$).
Since $|a-a_k| \to 0$ and $|b-b_k| \to 0$ for all $\omega'$,
and $f_k \to f$ uniformly in $\overline{\mathbb D}$,
we conclude that $|\tilde a - \tilde a_k| \to 0$
and $|\tilde b - \tilde b_k| \to 0$ for all $\omega'$.

Later we will prove a ``continuity'' property of $SLE_6$
(Lemma~\ref{continuity}) that allows us to conclude that,
under these conditions, $\tilde\gamma_k$ converges in
distribution to $\tilde\gamma$ in the uniform metric~(\ref{distance})
on continuous curves.
Once again, this implies the existence of coupled versions of
$\tilde\gamma_k$ and $\tilde\gamma$ on some probability space
$(\Omega', {\cal B}', {\mathbb P}')$ such that
$\text{d}(\tilde \gamma(\omega'),\tilde \gamma_k(\omega')) \to 0$,
for all $\omega'$ as $k \to \infty$.
Therefore, thanks to the convergence of $f_k$ to $f$
uniformly in $\overline{\mathbb D}$,
$\text{d}(f(\tilde\gamma(\omega')),f_k(\tilde\gamma_k(\omega'))) \to 0$,
for all $\omega'$ as $k \to \infty$.
But since $f(\tilde\gamma_k)$ is distributed as $\gamma_{D_k,a_k,b_k}$
and $f(\tilde\gamma)$ is distributed as $\gamma_{D,a,b}$,
we conclude that, as $k \to \infty$, $\gamma_{D_k,a_k,b_k}$ converges
in distribution to $\gamma_{D,a,b}$ in the uniform metric~(\ref{distance})
on continuous curves.

We now note that (S) implies that, as $\delta \to 0$,
$\gamma^{\delta}_{D_k,a_k,b_k}$ converges in distribution to
$\gamma_{D_k,a_k,b_k}$ \emph{uniformly} in $k$, for $k$ large enough.
Therefore, as $k \to \infty$, $\gamma^{\delta_k}_{D_k,a_k,b_k}$
converges in distribution to $\gamma_{D,a,b}$, and the proof is
concluded. \fbox{}

\begin{lemma} \label{continuity}
Let ${\mathbb D} \subset {\mathbb C}$ be the unit disc, $a$
and $b$ two distinct points on its boundary, and $\gamma$ the
trace of chordal $SLE_6$ inside $\mathbb D$ from $a$ to $b$.
Let $\{ a_k \}$ and $\{ b_k \}$ be two sequences of points in
$\partial {\mathbb D}$ such that $a_k \to a$ and $b_k \to b$.
Then, as $k \to \infty$, the trace $\gamma_k$ of chordal $SLE_6$
inside $\mathbb D$ from $a_k$ to $b_k$ converges in distribution
to $\gamma$ in the uniform topology~(\ref{distance}) on continuous
curves.
\end{lemma}

\noindent {\bf Proof.} Let
$f_k(z)=e^{i \alpha_k} \frac{z-z_k}{1-\bar z_k z}$
be the (unique) linear fractional transformation that takes
the unit disc $\mathbb D$ onto itself, mapping $a$ to $a_k$,
$b$ to $b_k$, and a third point $c \in \partial {\mathbb D}$
distinct from $a$ and $b$ to itself.
$\alpha_k$ and $z_k$ depend continuously on $a_k$ and $b_k$.
As $k \to \infty$, since $a_k \to a$ and $b_k \to b$,
$f_k$ converges uniformly to the identity in $\overline{\mathbb D}$.

Using the conformal invariance of chordal $SLE_6$, we couple
$\gamma_k$ and $\gamma$ by writing $\gamma_k=f_k(\gamma)$.
The uniform convergence of $f_k$ to the identity implies
that $\text{d}(\gamma,\gamma_k) \to 0$ as $k \to \infty$,
which is enough to conclude that $\gamma_k$ converges to
$\gamma$ in distribution. \fbox{} \\

\noindent {\bf Proof of Theorem~\ref{thm-convergence}.}
Let us prove the second part of the theorem first.
We will do this for the original version of the discrete
construction, but essentially the same proof works for the
reorganized version we will describe below, as we will explain
later.
Suppose that at step $k$ of this discrete construction
an exploration process $\gamma^{\delta}_k$ is run inside
a domain $D^{\delta}_{k-1}$, and write
$D^{\delta}_{k-1} \setminus \Gamma(\gamma^{\delta}_k) =
\bigcup_j D^{\delta}_{k,j}$, where $\{ D^{\delta}_{k,j} \}$
are the maximal connected domains of unexplored hexagons
into which $D^{\delta}_{k-1}$ is split by removing the
set $\Gamma(\gamma^{\delta}_k)$ of hexagons explored by
$\gamma^{\delta}_k$.

\begin{figure}[!ht]
\begin{center}
\includegraphics[width=8cm]{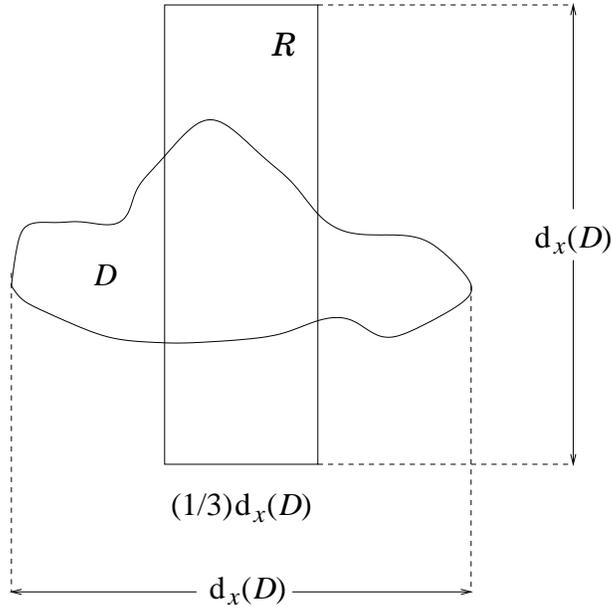}
\caption{Schematic drawing of a domain $D$ with
$\text{d}_x(D) \geq \text{d}_y(D)$ and the associated
rectangle $\cal R$.}
\label{fig3-thm1}
\end{center}
\end{figure}

Let $\text{d}_x(D^{\delta}_{k-1})$ and
$\text{d}_y(D^{\delta}_{k-1})$ be respectively the
maximal $x$- and $y$-distances between pairs of points
in $\partial D^{\delta}_{k-1}$.
Suppose, without loss of generality, that
$\text{d}_x(D^{\delta}_{k-1}) \geq \text{d}_y(D^{\delta}_{k-1})$,
and consider the rectangle $\cal R$ (see Figure~\ref{fig3-thm1})
whose vertical sides are aligned to the $y$-axis, have
length $\text{d}_x(D^{\delta}_{k-1})$, and are each
placed at $x$-distance
$\frac{1}{3} \, \text{d}_x(D^{\delta}_{k-1})$ from
points of $\partial D^{\delta}_{k-1}$ with minimal
or maximal $x$-coordinate in such a way that the
horizontal sides of $\cal R$ have length
$\frac{1}{3} \, \text{d}_x(D^{\delta}_{k-1})$;
the bottom and top sides of $\cal R$ are placed in
such a way that they are at equal $y$-distance from
the points of $\partial D^{\delta}_{k-1}$ with minimal
or maximal $y$-coordinate, respectively.

It follows from the Russo-Seymour-Welsh lemma~\cite{russo,sewe}
(see also~\cite{kesten,grimmett}) that the probability to
have two vertical $\cal T$-crossings of $\cal R$ of different
colors is bounded away from zero by a positive constant $p_0$
that does not depend on $\delta$ (for $\delta$ small enough).
If that happens, then
$\max_j \text{d}_x(D^{\delta}_{k,j}) \leq \frac{2}{3} \text{d}_x(D^{\delta}_{k-1})$.
The same argument of course applies to the maximal
$y$-distance when
$\text{d}_y(D^{\delta}_{k-1}) \geq \text{d}_x(D^{\delta}_{k-1})$.
We can summarize the above observation in the following lemma.
\begin{lemma} \label{smaller}
Suppose that at step $k$ of the full discrete construction
an exploration process $\gamma^{\delta}_k$ is run inside
a domain $D^{\delta}_{k-1}$.
If $\text{\emph{d}}_x(D^{\delta}_{k-1}) \geq \text{\emph{d}}_y(D^{\delta}_{k-1})$,
then for $\delta$ small enough (i.e., $\delta \leq C \, \text{\emph{d}}_x(D^{\delta}_{k-1})$
for some constant $C$),
$\max_j \text{\emph{d}}_x(D^{\delta}_{k,j}) \leq \frac{2}{3} \text{\emph{d}}_x(D^{\delta}_{k-1})$
with probability at least $p_0$ independent of $\delta$.
The same holds for the maximal $y$-distances when
$\text{\emph{d}}_y(D^{\delta}_{k-1}) \geq \text{\emph{d}}_x(D^{\delta}_{k-1})$.
\end{lemma}

%

Here is another lemma that will be useful later on.
(For an example of the phenomenon described in the lemma, see
Figure~\ref{fig2-sec4}, and assume that the unexplored hexagons
there are all blue; then the s-boundary of the small domain
made of a single blue hexagon and that of the blue domain to
the northeast share exactly two adjacent yellow hexagons.)
\begin{lemma} \label{daughters}
Two ``daughter" subdomains, $D^{\delta}_{k,j}$ and $D^{\delta}_{k,j'}$,
either have disjoint s-boun\-da\-ries, or else their common s-boundary
consists of exactly two adjacent hexagons (of the same color) where
the exploration path $\gamma^{\delta}_k$ came within $2$ hexagons of
touching itself just when completing the s-boundary of one of the
two subdomains.
\end{lemma}

\noindent {\bf Proof.} Suppose that the two daughter subdomains
have s-boundaries $\Delta D^{\delta}_{k,j}$ and $\Delta D^{\delta}_{k,j'}$
that are not disjoint and let $S = \{ \xi_1,\ldots,\xi_i \}$ be the set
of (sites of $\cal T$ that are the centers of the) hexagons that belong
to both s-boundaries.
$S$ can be partitioned into subsets consisting of single hexagons
that are not adjacent to any another hexagon in $S$ and groups of
hexagons that form simple $\cal T$-paths (because the s-boundaries
of the two subdomains are simple $\cal T$-loops).
Let $\{ \xi_l,\ldots,\xi_m \}$ be such a subset of hexagons of $S$
that form a simple $\cal T$-path $\pi_0=(\xi_l,\ldots,\xi_m)$.
Then there is a $\cal T$-path $\pi_1$ of hexagons in $\Delta D^{\delta}_{k,j}$
that goes from $\xi_l$ to $\xi_m$ without using any other hexagon of $\pi_0$
and a different $\cal T$-path $\pi_2$ in $\Delta D^{\delta}_{k,j'}$ that
goes from $\xi_m$ to $\xi_l$ without using any other hexagon of $\pi_0$.
But then, all the hexagons in $\pi_0$ other than $\xi_l$ and $\xi_m$ are
``surrounded" by $\pi_1 \cup \pi_2$ and therefore cannot have been explored
by the exploration process that produced $D^{\delta}_{k,j}$ and $D^{\delta}_{k,j'}$,
and cannot belong to $\Delta D^{\delta}_{k,j}$ or $\Delta D^{\delta}_{k,j'}$,
leading to a contradiction, unless $\pi_0 = (\xi_l,\xi_m)$.
Similar arguments lead to a contradiction if $S$ is partitioned into more
than one subset.

If $\xi_i \in S$ is not adjacent to any other hexagon in $S$, then it
is adjacent to two other hexagons of $\Delta D^{\delta}_{k,j}$ and two
hexagons of $\Delta D^{\delta}_{k,j'}$.
Since $\xi_i$ has only six neighbors and neither the two hexagons of
$\Delta D^{\delta}_{k,j}$ adjacent to $\xi_i$ nor those of $\Delta D^{\delta}_{k,j'}$
can be adjacent to each other, each hexagon of $\Delta D^{\delta}_{k,j}$
is adjacent to one of $\Delta D^{\delta}_{k,j'}$.
But then, as before, $\xi_i$ is ``surrounded" by
$\{ \Delta D^{\delta}_{k,j} \cup \Delta D^{\delta}_{k,j'} \} \setminus \xi_i$
and therefore cannot have been explored by the exploration process that produced
$D^{\delta}_{k,j}$ and $D^{\delta}_{k,j'}$, and cannot belong to $\Delta D^{\delta}_{k,j}$
or $\Delta D^{\delta}_{k,j'}$, leading once again to a contradiction.
The proof is now complete, since the only case remaining is the one where
$S$ consists of a single pair of adjacent hexagons as stated in the lemma. \fbox{} \\

With these lemmas, we can now proceed with the proof of the second
part of the theorem.
Lemma~\ref{smaller} tells us that large domains are ``chopped" with
bounded away from zero probability ($\geq p_0 > 0$), but we need to
keep track of domains of diameter larger than $\varepsilon$ in such
a way as to avoid ``double counting" as the lattice construction proceeds.
More accurately, we will keep track of domains $\tilde D^{\delta}$ having
$\text{d}_m(\tilde D^{\delta}) \geq \frac{1}{\sqrt{2}} \, \varepsilon$,
since only these can have diameter larger than $\varepsilon$.
To do so, we will associate with each domain $\tilde D^{\delta}$ having
$\text{d}_m(\tilde D^{\delta}) \geq \frac{1}{\sqrt{2}} \, \varepsilon$
that we encounter as we do the lattice construction a non-negative
integer label.
The first domain is $D_0^{\delta} = {\mathbb D}^{\delta}$ (see the beginning
of Section~\ref{full}) and this gets label $1$.
After each exploration process in a domain $\tilde D^{\delta}$ with
$\text{d}_m(\tilde D^{\delta}) \geq \frac{1}{\sqrt{2}} \, \varepsilon$,
if the number $\tilde m$ of ``daughter" subdomains $\tilde D^{\delta}_j$ with
$\text{d}_m(\tilde D_j^{\delta}) \geq \frac{1}{\sqrt{2}} \, \varepsilon$
is $0$, then the label of $\tilde D^{\delta}$ is no longer used, if instead
$\tilde m \geq 1$, then one of these $\tilde m$ subdomains (chosen by any
procedure -- e.g., the one with the highest priority for further exploration)
is assigned the \emph{same} label as $\tilde D^{\delta}$ and the rest are
assigned the next $\tilde m - 1$ integers that have never before been used
as labels.
Note that once all domains have $\text{d}_m < \frac{1}{\sqrt{2}} \, \varepsilon$,
there are no more labelled domains.

%

\begin{lemma} \label{labels}
Let $M^{\delta}_{\varepsilon}$ denote the total number of labels
used in the above procedure; then for any fixed $\varepsilon>0$,
$M^{\delta}_{\varepsilon}$ is bounded in probability as $\delta \to 0$;
i.e.,
$\lim_{M \to \infty} \limsup_{\delta \to 0}
{\mathbb P}(M^{\delta}_{\varepsilon} > M) = 0$.
\end{lemma}

\noindent {\bf Proof.} Except for $D^{\delta}_0$, every domain comes with
(at least) a ``physically correct" monochromatic ``half-boundary" (notice
that we are considering s-boundaries and that a half-boundary coming from
the ``artificially colored" boundary of $D^{\delta}_0$ is not considered
a physically correct monochromatic half-boundary).
Let us assume, without loss of generality, that $M^{\delta}_{\varepsilon}>1$.
If we associate with each label the ``last" (in terms of steps of the discrete
construction) domain which used that label (its daughter subdomains all had
$\text{d}_m < \frac{1}{\sqrt{2}} \, \varepsilon$), then we claim that it
follows from Lemma~\ref{daughters} that (with high probability) any two such
last domains that are labelled have disjoint s-boundaries.
This is a consequence of the fact that the two domains are subdomains of two
``ancestors" that are distinct daughter subdomains of the \emph{same} domain
(possibly $D^{\delta}_0$) and whose s-boundaries are therefore (by Lemma~\ref{daughters})
either disjoint or else overlap at a pair of hexagons where an exploration path
had a close encounter of distance two hexagons with itself.
But since we are dealing only with macroscopic domains (of diameter at least
order $\varepsilon$), such a close encounter would imply, like in Lemmas~\ref{sub-conv}
and~\ref{boundaries}, the existence of six crossings, not all of the same color,
of an annulus whose outer radius can be kept fixed while the inner radius is
sent to zero together with $\delta$.
The probability of such an event goes to zero as $\delta \to 0$ and hence
the unit disc $\mathbb D$ contains, with high probability, at least
$M^{\delta}_{\varepsilon}$ disjoint monochromatic $\cal T$-paths of diameter
at least $\frac{1}{\sqrt{2}} \, \varepsilon$, corresponding to the physically
correct half-boundaries of the $M^{\delta}_{\varepsilon}$ labelled domains.

Now take the collection of squares $s_j$ of side length
$\varepsilon'>0$ centered at the sites $c_j$ of a scaled
square lattice $\varepsilon' {\mathbb Z}^2$ of mesh size
$\varepsilon'$, and let $N(\varepsilon')$ be the number
of squares of side $\varepsilon'$ needed to cover the
unit disc.
Let $\varepsilon'<\varepsilon/2$ and consider the
event $\{ M^{\delta}_{\varepsilon} \geq 6 \, N(\varepsilon') \}$, which
implies that, with high probability, the unit disc contains at least
$6 \, N(\varepsilon')$ disjoint monochromatic $\cal T$-paths of diameter at
least $\frac{1}{\sqrt{2}} \, \varepsilon$ and that, for at least one $j=j_0$,
the square $s_{j_0}$ intersects at least six disjoint monochromatic
$\cal T$-paths of diameter larger that $\frac{1}{\sqrt{2}} \, \varepsilon$,
so that the ``annulus" $B(c_{j_0},\frac{1}{2 \sqrt{2}} \, \varepsilon) \setminus s_{j_0}$
is crossed by at least six disjoint monochromatic $\cal T$-paths contained
inside the unit disc.

If all these $\cal T$-paths crossing
$B(c_{j_0},\frac{1}{2 \sqrt{2}} \, \varepsilon) \setminus s_{j_0}$
have the same color, say blue, then since they are portions of boundaries
of domains discovered by exploration processes, they are ``shadowed" by
exploration paths and therefore between at least one pair of blue $\cal T$-paths,
there is at least one yellow $\cal T$-path crossing
$B(c_{j_0},\frac{1}{2 \sqrt{2}} \, \varepsilon) \setminus s_{j_0}$.
Therefore, whether the original monochromatic $\cal T$-paths are
all of the same color or not,
$B(c_{j_0},\frac{1}{2 \sqrt{2}} \, \varepsilon) \setminus s_{j_0}$
is crossed by at least six disjoint monochromatic $\cal T$-paths
not all of the same color contained in the unit disc.
%
%
Let $g(\varepsilon,\varepsilon')$ denote the $\limsup$ as $\delta \to 0$
of the probability that such an event happens anywhere inside the unit disc.
We have shown that the event $\{ M^{\delta}_{\varepsilon} \geq 6 \, N(\varepsilon') \}$
implies a ``six-arms" event unless not all labelled domains have disjoint s-boundaries.
But the latter also implies a ``six-arms" event, as discussed before; therefore
\begin{equation}
\limsup_{\delta \to 0} {\mathbb P}(M^{\delta}_{\varepsilon} \geq 6 \, N(\varepsilon'))
\leq 2 \, g(\varepsilon,\varepsilon').
\end{equation}
Since $B(c_{j_0},\frac{1}{2 \sqrt{2}} \, \varepsilon) \setminus B(c_{j_0},\frac{1}{\sqrt{2}} \, \varepsilon')
\subset B(c_{j_0},\frac{1}{2 \sqrt{2}} \, \varepsilon) \setminus s_{j_0}$,
bounds in~\cite{ksz} imply that, for $\varepsilon$ fixed,
$g(\varepsilon,\varepsilon') \to 0$ as $\varepsilon' \to 0$,
which shows that
\begin{equation}
\lim_{M \to \infty} \limsup_{\delta \to 0} {\mathbb P}(M^{\delta}_{\varepsilon} > M) = 0
\end{equation}
and concludes the proof of the lemma. \fbox{} \\

Now, let $N^{\delta}_i$ denote the number of distinct domains that had
label $i$ (this is equal to the number of steps that label $i$ survived).
Let us also define $H(\varepsilon)$ to be the smallest integer $h \geq 1$
such that $(\frac{2}{3})^h<\frac{1}{\sqrt{2}} \, \varepsilon$ and $G_h$ to
be the random variable corresponding to how many Bernoulli trials (with
probability $p_0$ of success) it takes to have $h$ successes.
Then, we may apply (sequentially) Lemma~\ref{smaller} to conclude that
for any $i$
\begin{equation}
{\mathbb P}(N^{\delta}_i \geq k+1) \leq {\mathbb P}(G_{H(\varepsilon)} + G'_{H(\varepsilon)} \geq k),
\end{equation}
where $G'_h$ is an independent copy of $G_h$.

Now let $\tilde N_1(\varepsilon),\tilde N_2(\varepsilon),\ldots$ be i.i.d.
random variables equidistributed with $G_{H(\varepsilon)} + G'_{H(\varepsilon)}$.
Let $\tilde K_{\delta}(\varepsilon)$ be the number of steps needed so that
all domains left to explore have $\text{d}_m<\frac{1}{\sqrt{2}} \, \varepsilon$.
Then, for any positive integer $M$,
\begin{equation}
{\mathbb P}(\tilde K_{\delta}(\varepsilon)>C)
\leq {\mathbb P}(M^{\delta}_{\varepsilon} \geq M+1)
+{\mathbb P}(\tilde N_1(\varepsilon)+\ldots+\tilde N_M(\varepsilon) \geq C).
\end{equation}
Notice that, for fixed $M$,
${\mathbb P}(\tilde N_1(\varepsilon)+\ldots+\tilde N_M(\varepsilon) \geq C) \to 0$
as $C \to \infty$.
Moreover, for any $\hat\varepsilon>0$, by Lemma~\ref{labels},
we can choose $M_0=M_0(\hat\varepsilon)$ large enough so that
$\limsup_{\delta \to 0} {\mathbb P}(M^{\delta}_{\varepsilon} > M_0) < \hat\varepsilon$.
So, for any $\hat\varepsilon>0$, it follows that
\begin{equation}
\limsup_{C \to \infty} \limsup_{\delta \to 0} {\mathbb P}(\tilde K_{\delta}(\varepsilon)>C)
< \hat\varepsilon,
\end{equation}
which implies that
\begin{equation}
\lim_{C \to \infty} \limsup_{\delta \to 0} {\mathbb P}(\tilde K_{\delta}(\varepsilon)>C)
= 0.
\end{equation}

To conclude this part of the proof, notice that the discrete construction
cannot ``skip'' a contour and move on to explore its interior, so that all the
contours with diameter larger than $\varepsilon$ must have been found by step
$k$ if all the domains present at that step have diameter smaller than $\varepsilon$.
Therefore, $K_{\delta}(\varepsilon) \leq \tilde K_{\delta}(\varepsilon)$, which
shows that $K_{\delta}(\varepsilon)$ is bounded in probability as $\delta \to 0$. \\

For the first part of the theorem, we need to prove, for
any fixed $k \in {\mathbb N}$, joint convergence in distribution
of the first $k$ steps of a suitably reorganized discrete
construction to the first $k$ steps of the continuum one.
Later we will explain why this reorganized construction has
the same scaling limit as the one defined in Section~\ref{full}.
For each $k$, the first $k$ steps of the reorganized discrete
construction will be coupled to the first $k$ steps of the continuum
one with suitable couplings in order to obtain the convergence in
distribution of those steps of the discrete construction to the analogous
steps of the continuum one; the proof will proceed by induction in $k$.
We will explain how to reorganize the discrete construction as we go along;
in order to explain the idea of the proof, we will consider first the
cases $k=1$, $2$ and $3$, and then extend to all $k>3$. \\


\noindent $k=1$. The first step of the continuum construction
consists of an $SLE_6$ $\gamma_1$ from $-i$ to $i$ inside
$\mathbb D$.
Correspondingly, the first step of the discrete construction
consists of an exploration path $\gamma^{\delta}_1$ inside
${\mathbb D}^{\delta}$ from the e-vertex closest to $-i$ to
the e-vertex closest to $i$.
The convergence in distribution of $\gamma^{\delta}_1$ to
$\gamma_1$ is covered by statement (S). \\

\noindent $k=2$.
The convergence in distribution of the percolation exploration
path to chordal $SLE_6$ implies that we can couple $\gamma^{\delta}_1$
and $\gamma_1$ generating them as random variables on some probability
space $(\Omega',{\cal B}',{\mathbb P}')$
such that
$\text{d}(\gamma_1(\omega'),\gamma^{\delta}_1(\omega')) \to 0$
for all $\omega'$ as $k \to \infty$ (see, for example,
Corollary~1 of~\cite{billingsley1}).

Now, let $D_1$ be the domain generated by $\gamma_1$ that
is chosen for the second step of the continuum construction,
and let $c_1 \in {\cal P}$ be the highest ranking
point of $\cal P$ contained in $D_1$.
For $\delta$ small enough, $c_1$ is also contained in
${\mathbb D}^{\delta}$; let $D_1^{\delta} = D_1^{\delta}(c_1)$
be the unique connected component of the set
${\mathbb D}^{\delta} \setminus \Gamma(\gamma_1^{\delta})$
containing $c_1$ (this is well-defined with probability
close to $1$ for small $\delta$); $D_1^{\delta}$ is the domain
where the second exploration process is to be carried out.
>From the proof of Lemma~\ref{boundaries}, we know that
the boundaries $\partial D_1^{\delta}$ and $\partial D_1$
of the domains $D_1^{\delta}$ and $D_1$ produced respectively
by the path $\gamma_1^{\delta}$ and $\gamma_1$ are close
with probability close to one for $\delta$ small enough.


For the next step of the discrete construction, we choose
the two e-vertices $x_1$ and $y_1$ in $\partial D_1^{\delta}$
that are closest to the points $a_1$ and $b_1$ of $\partial D_1$
selected for the coupled continuum construction (if the choice
is not unique, we can select the e-vertices with any rule to
break the tie) and call $\gamma^{\delta}_2$ the percolation
exploration path inside $D^{\delta}_1$ from $x_1$ to $y_1$.
It follows from ~\cite{ab} that
$\{ \gamma_1^{\delta},\partial D^{\delta}_1,\gamma^{\delta}_2 \}$
converge jointly in distribution along some subsequence
to some limit
$\{ \tilde\gamma_1,\partial \tilde D_1,\tilde\gamma_2 \}$.
We already know that $\tilde\gamma_1$ is distributed like $\gamma_1$
and we can deduce from the joint convergence in distribution of
$(\gamma^{\delta}_1,\partial D_1^{\delta})$ to $(\gamma_1,\partial D_1)$
(Lemma~\ref{boundaries}), that $\partial \tilde D_1$ is distributed
like $\partial D_1$.
Therefore, if we call $\gamma_2$ the $SLE_6$ path inside $D_1$ from
$a_1$ to $b_1$, Lemma~\ref{strong-smirnov}
implies that $\tilde\gamma_2$ is distributed like $\gamma_2$ and indeed
that, as $\delta \to 0$,
$\{ \gamma^{\delta}_1,\partial D^{\delta}_1,\gamma^{\delta}_2 \}$ converge
jointly in distribution to $\{ \gamma_1,\partial D_1,\gamma_2 \}$. \\

\noindent $k = 3$. So far, we have proved the convergence
in distribution of the (paths and boundaries produced in the)
first two steps of the discrete construction to the (paths and
boundaries produced in the) first two steps of the discrete
construction.
The third step of the continuum construction consists of an $SLE_6$
path $\gamma_3$ from $a_2 \in \partial D_2$ to $b_2 \in \partial D_2$,
inside the domain $D_2$ with highest priority after the second step
has been completed.
Let $c_2 \in {\cal P}$ be the highest ranking point of $\cal P$ contained
in $D_2$, $D_2^{\delta}$ the domain of the discrete construction containing
$c_2$ after the second step of the discrete construction has been completed
(this is well defined with probability close to $1$ for small $\delta$),
and choose the two e-vertices $x_2$ and $y_2$ in $\partial D_2^{\delta}$
that are closest to the points $a_2$ and $b_2$ of $\partial D_2$ selected
for the coupled continuum construction (if the choice is not unique, we
can select the e-vertices with any rule to break the tie).
The third step of the discrete construction consists of an exploration
path $\gamma_3^{\delta}$ from $x_2$ to $y_2$ inside $D_2^{\delta}$.

It follows from ~\cite{ab} that
$\{ \gamma_1^{\delta},\partial D^{\delta}_1,\gamma^{\delta}_2,\partial D_2^{\delta},\gamma_3^{\delta} \}$
converge jointly in distribution along some subsequence
to some limit
$\{ \tilde\gamma_1,\partial\tilde D_1,\tilde\gamma_2,\partial\tilde D_2,\tilde\gamma_3 \}$.
We already know that $\tilde\gamma_1$ is distributed like $\gamma_1$,
$\partial\tilde D_1$ like $\partial D_1$ and $\tilde\gamma_2$ like $\gamma_2$,
and we would like to apply Lemma~\ref{strong-smirnov} to conclude
that $\tilde\gamma_3$ is distributed like $\gamma_3$
and indeed that, as $\delta \to 0$,
$(\gamma_1^{\delta},\partial D^{\delta}_1,\gamma^{\delta}_2,\partial D_2^{\delta},\gamma_3^{\delta})$
converges in distribution to $(\gamma_1,\partial D_1,\gamma_2,\partial D_2,\gamma_3)$.
In order to do so, we have to first show that $\partial\tilde D_2$
is distributed like $\partial D_2$.
If $D^{\delta}_2$ is a subset of ${\mathbb D}^{\delta} \setminus \Gamma(\gamma^{\delta}_1)$,
this follows from Lemma~\ref{boundaries}, as in the previous case,
but if the s-boundary of $D^{\delta}_2$ contains hexagons of
$\Gamma(\gamma^{\delta}_2)$, then we cannot use Lemma~\ref{boundaries}
directly, although the proof of the lemma can be easily adapted to
the present case, as we now explain.

\begin{figure}[!ht]
\begin{center}
\includegraphics[width=8cm]{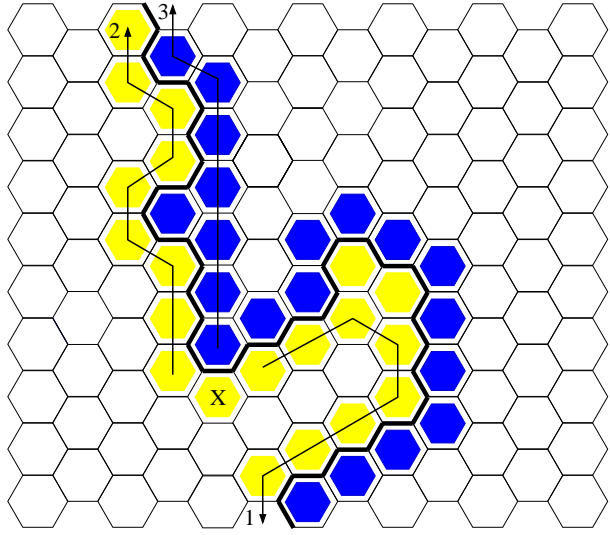}
\caption{Hexagon X, in the s-boundary of the domain $D_j^{\delta}$
to the left of the exploration path indicated by a heavy line,
has three neighbors that are the starting points of two
disjoint yellow $\cal T$-paths (denoted 1 and 2) belonging to
the s-boundary of $D_j^{\delta}$ and one blue $\cal T$-path
(denoted 3) belonging to the adjacent percolation cluster.}
\label{fig2-thm1}
\end{center}
\end{figure}

Indeed, the only difference is in the proof of claim (C) and is
due to the fact that, when dealing with a domain of type~1 or 2,
we cannot use the bound on the probability of three disjoint
crossings of a semi-annulus because the domains we are dealing
with may not be convex (like the unit disc).
On the other hand, the discrete domains like $D^{\delta}_1$ and
$D^{\delta}_2$ where we have to run exploration processes at
various steps of the discrete construction are themselves generated
by previous exploration processes, so that any hexagon of the
s-boundary of such a domain has three adjacent hexagons which are
the starting points of three disjoint $\cal T$-paths (two of
one color and one of the other).
Two of these $\cal T$-paths belong to the s-boundary
of the domain, while the third belongs to the adjacent
percolation cluster (see Figure~\ref{fig2-thm1}).
This allows us to use the bound on the probability of
\emph{six} disjoint crossings of an annulus.

To see this, let $\pi_1, \pi_2$ be the $\cal T$-paths
contained in the s-boundary of the discrete domain (i.e.,
$D^{\delta}_1$ in the present context) and $\pi_3$ the
$\cal T$-path belonging to the adjacent cluster, all starting
from hexagons adjacent to some hexagon $\xi$ (centered at $u$)
in the s-boundary of $D^{\delta}_1$.
For $0< \varepsilon' < \varepsilon$ and $\delta$ small
enough, let ${\cal A}_u(\varepsilon,\varepsilon')$ be the
event that the exploration path $\gamma_2^{\delta}$
enters the ball $B(u,\varepsilon')$ without touching
$\partial D^{\delta}_1$ inside the larger ball $B(u,\varepsilon)$.
${\cal A}_u(\varepsilon,\varepsilon')$ implies having (at least)
three disjoint $\cal T$-paths (two of one color and one of the other),
$\pi_4, \pi_5$ and $\pi_6$, contained in $D^{\delta}_1$ and crossing
the annulus $B(u,\varepsilon) \setminus B(u,\varepsilon')$, with
$\pi_4, \pi_5$ and $\pi_6$ disjoint from $\pi_1, \pi_2$ and $\pi_3$.
Hence, ${\cal A}_u(\varepsilon,\varepsilon')$ implies the event that
there are (at least) six disjoint crossings (not all of the same color)
of the annulus $B(u,\varepsilon) \setminus B(u,\varepsilon')$.

Once claim (C) is proved, the rest of the proof of
Lemma~\ref{boundaries} applies to the present case.
Therefore, we have convergence in distribution of
$\partial D_2^{\delta}$ to $\partial D_2$, which allows
us to use Lemma~\ref{strong-smirnov} and conclude that
$(\gamma_1^{\delta},\gamma_2^{\delta},\gamma_3^{\delta})$
converges in distribution to $(\gamma_1,\gamma_2,\gamma_3)$. \\

\noindent $k > 3$.
We proceed by induction in $k$, iterating the steps explained
above; there are no new difficulties; all steps for $k \geq 4$
are analogous to the case $k = 3$. \\

To conclude the proof of the theorem, we need to show
that the scaling limit of the original full discrete
construction defined in Section~\ref{full} is the same
as that of the reorganized one just used in the proof
of the first part of the theorem.
In order to do so, we can couple the two constructions
by using the same percolation configuration for both,
so that the two constructions have at their disposal
the same set of loops to discover.
We proved above that the original discrete construction
finds \emph{all} the ``macroscopic" loops, so we have to
show that this is true also for the reorganized version
of the discrete construction.
This is what we will do next, using essentially the same
arguments as those employed for the original discrete
construction; we present these arguments for the sake
of completeness since there are some changes.

Consider the reorganized discrete construction described
above, where the starting and ending points of the exploration
processes at each step are chosen to be close to those of
the corresponding (coupled) continuum construction.
Suppose that at step $k$ of this discrete construction
an exploration process $\gamma^{\delta}_k$ is run inside
a domain $D^{\delta}_{k-1}$, and write
$D^{\delta}_{k-1} \setminus \Gamma(\gamma^{\delta}_k) =
\bigcup_j D^{\delta}_{k,j}$, where $\{ D^{\delta}_{k,j} \}$
are the connected domains into which $D^{\delta}_{k-1}$
is split by the set $\Gamma(\gamma^{\delta}_k)$ of hexagons
explored by $\gamma^{\delta}_k$.

Let $\text{d}_x(D_{k-1})$ (resp., $\text{d}_x(D^{\delta}_{k-1})$)
and $\text{d}_y(D_{k-1})$ (resp., $\text{d}_y(D^{\delta}_{k-1})$)
be respectively the maximal $x$- and $y$-distance between pairs
of points in $\partial D_{k-1}$ (resp., $\partial D^{\delta}_{k-1}$).
%
%
%
If
$\text{d}_x(D^{\delta}_{k-1}) \geq \text{d}_y(D^{\delta}_{k-1})$
and the e-vertices on $\partial D^{\delta}_{k-1}$
are chosen to be closest to two points of $\partial D_{k-1}$
with maximal $x$-distance, then the same construction and
argument spelled out earlier in the first part of the proof
(corresponding to the second part of the theorem) show that
$\max_j \text{d}_x(D^{\delta}_{k,j}) \leq \frac{2}{3} \text{d}_x(D^{\delta}_{k-1})$
with bounded away from zero probability.

%

If the e-vertices on $\partial D^{\delta}_{k-1}$ are
chosen to be closest to two points of $\partial D_{k-1}$
with maximal $x$-distance but
$\text{d}_x(D^{\delta}_{k-1}) \leq \text{d}_y(D^{\delta}_{k-1})$,
then consider the rectangle ${\cal R}'$ whose vertical sides
are aligned to the $y$-axis, have length $\text{d}_y(D^{\delta}_{k-1})$,
and are each placed at the same $x$-distance from the
points of $\partial D^{\delta}_{k-1}$ with minimal
or maximal $x$-coordinate in such a way that the
horizontal sides of ${\cal R}'$ have length
$\frac{1}{3} \, \text{d}_y(D^{\delta}_{k-1})$;
the bottom and top sides of ${\cal R}'$ are placed in
such a way that they touch the points of
$\partial D^{\delta}_{k-1}$ with minimal or maximal
$y$-coordinate, respectively.
Notice that, because of the coupling between the continuum
and discrete constructions, for any $\tilde\varepsilon >0$, for
$k$ large enough, $|\text{d}_x(D_{k-1}^{\delta}) - \text{d}_x(D_{k-1})| \leq \tilde\varepsilon$
and $|\text{d}_y(D_{k-1}^{\delta}) - \text{d}_y(D_{k-1})| \leq \tilde\varepsilon$.
Since in the case under consideration we have
$\text{d}_y(D_{k-1}^{\delta}) \geq \text{d}_x(D_{k-1}^{\delta})$ and
$\text{d}_x(D_{k-1}) \geq \text{d}_y(D_{k-1})$, for $\delta$ large enough,
we must also have
$|\text{d}_y(D_{k-1}^{\delta}) - \text{d}_x(D_{k-1}^{\delta})| \leq 2 \, \tilde\varepsilon$.
Once again, it follows from the Russo-Seymour-Welsh lemma
that the probability to have two vertical $\cal T$-crossings
of ${\cal R}'$ of different colors is bounded away from zero by
a positive constant that does not depend on $\delta$ (for $\delta$
small enough).
If that happens, then
$\max_j \text{d}_x(D^{\delta}_{k,j}) \leq
\frac{2}{3} \, \text{d}_x(D^{\delta}_{k-1}) + \frac{1}{3} \, \tilde\varepsilon$.

All other cases are handled in the same way, implying that
the maximal $x$- and $y$-distances of domains that appear
in the discrete construction have a positive probability
(bounded away from zero) to decrease by (approximately) a
factor $2/3$ at each step of the discrete construction in
which an exploration process is run in that domain.


With this result at our disposal, the rest of the proof, that
for any $\varepsilon>0$ the number of steps needed to find
all the loops of diameter larger than $\varepsilon$ is bounded
in probability as $\delta \to 0$ (which implies that \emph{all}
the ``macroscopic" loops are discovered), proceeds exactly like
for the original discrete construction. \fbox{} \\

\noindent {\bf Proof of Theorem~\ref{thm-therm-lim}.}
First of all, we want to show that $P^D \equiv \hat I_D P_R$
does not depend on $R$, provided $D$ is strictly contained
in ${\mathbb D}_R$ and
$\partial D \cap \partial {\mathbb D}_R = \emptyset$.
In order to do this, we assume that the above
conditions are satisfied for the pair $D,R$ and
show that $\hat I_D P_R = \hat I_D P_{R'}$ for all $R' > R$.

Take two copies of the scaled hexagonal lattice,
$\delta{\cal H}$ and $\delta{\cal H}'$, their
dual lattices $\delta{\cal T}$ and $\delta{\cal T}'$,
and two percolation configurations,
$\sigma_{{\mathbb D}_R}$ and ${\sigma'}_{{\mathbb D}_{R'}}$,
both with blue boundary conditions and coupled in
such a way that
$\sigma_{{\mathbb D}_R} = \sigma'_{{\mathbb D}_R}$.
The laws of the boundaries of $\sigma$ and $\sigma'$
are also coupled, in such a way that the boundaries
or portions of boundaries contained inside $D$ are
identical for all small enough $\delta$.
Therefore, letting $\delta \to 0$ and using the
convergence of the percolation boundaries inside
${\mathbb D}_R$ and ${\mathbb D}_{R'}$ to the
continuum nonsimple loop processes $P_R$ and $P_{R'}$
respectively, we conclude that $\hat I_D P_R = \hat I_D P_{R'}$.

>From what we have just proved, it follows that the probability
measures $P^{{\mathbb D}_R}$ on $(\Omega_R,{\cal B}_R)$, for
$R \in {\mathbb R}_+$, satisfy the consistency conditions
$P^{{\mathbb D}_{R_1}} = \hat I_{{\mathbb D}_{R_1}} P^{{\mathbb D}_{R_2}}$
for all $R_1 \leq R_2$.
Since $\Omega_R$, $\Omega$ are complete separable metric
spaces, the measurable spaces $(\Omega_R,{\cal B}_R)$,
$(\Omega,{\cal B})$ are standard Borel spaces and so we
can apply Kolmogorov's extension theorem (see, for example,
\cite{durrett}) and conclude that there exists a unique
probability measure on $(\Omega,{\cal B})$ with
$P^{{\mathbb D}_R} = \hat I_{{\mathbb D}_R} P$ for all
$R \in {\mathbb R}_+$.
It follows that, for $R'>R$ and all $D$ strictly contained
in ${\mathbb D}_R$ and such that
$\partial D \cap \partial {\mathbb D}_R = \emptyset$,
$\hat I_D P_R = P^D = \hat I_D P_{R'} =
\hat I_D \hat I_{{\mathbb D}_R} P_{R'} =
\hat I_D P^{{\mathbb D}_R} =
\hat I_D \hat I_{{\mathbb D}_R} P = \hat I_D P$, which
concludes the proof. \fbox{} \\

\noindent {\bf Proof of Theorems~\ref{main-result1} and~\ref{main-result2}.}
These are immediate consequences of
Theorems~\ref{thm-convergence} and~\ref{thm-therm-lim},
where the full scaling limit is intended in the topology
induced by~(\ref{hausdorff-D}). \fbox{} \\

\noindent {\bf Proof of Theorem~\ref{features}.}
{\it 1.} The fact that the Continuum Nonsimple Loop process
is a random collection of noncrossing continuous loops is a
direct consequence of its definition.
The fact that the loops touch themselves is a
consequence of their being constructed out of $SLE_6$,
while the fact that they touch each other follows
from the observation that a chordal $SLE_6$ path
$\gamma_{D,a,b}$ touches $\partial D$ with probability
one.
Therefore, each new loop in the continuum construction
touches one or more previous ones (many times).

The nonexistence of triple points follows directly from
Lemma~5 of~\cite{ksz} on the number of crossings of an
annulus, combined with Theorem~\ref{main-result1}, which
allows to transport discrete results to the continuum case.
In fact, a triple point would imply, for discrete percolation,
at least six crossings (not all of the same color) of an
annulus whose ratio of inner to outer radius goes to zero
in the scaling limit, leading to a contradiction.

{\it 2.} This follows from straightforward
Russo-Seymour-Welsh type arguments for percolation
(for more details, see, for example, Lemma~3
of~\cite{ksz}), combined with Theorem~\ref{main-result1}.

{\it 3.} 
Combining Russo-Seymour-Welsh type arguments for
percolation (see, for example, Lemma~3 of~\cite{ksz})
with Theorem~\ref{main-result1}, we know that $P$-a.s.
there exists a (random) $R^* = R^*(R)$, with $R^* < \infty$,
such that ${\mathbb D}_R$ is surrounded by a continuum
nonsimple loop contained in ${\mathbb D}_{R^*}$.
>From (the proof of) Theorem~\ref{thm-therm-lim}, we also
know that
$\hat I_{{\mathbb  D}_{R''}} P = P^{{\mathbb D}_{R''}} =
I_{{\mathbb D}_{R''}} P_{R'}$ for all $R' > R''$.
This implies that by taking $R'$ large enough
and performing the continuum construction inside ${\mathbb D}_{R'}$,
we have a positive probability of generating a loop $\lambda$
contained in the annulus ${\mathbb D}_{R''} \setminus {\mathbb D}_R$,
with $R'>R''>R$.
If that is the case, all the loops contained inside
${\mathbb D}_R$ are connected, by construction, to the
loop $\lambda$ surrounding ${\mathbb D}_R$ by a finite
sequence (a ``path'') of loops (remember that in the
continuum construction each loop is generated by pasting
together portions of $SLE_6$ paths inside domains whose
boundaries are determined by previously formed loops or
excursions).
Therefore, any two loops contained inside ${\mathbb D}_R$
are connected to each other by a ``path" of loops.

Using again the fact that
$\hat I_{{\mathbb  D}_{R''}} P = P^{{\mathbb D}_{R''}} =
I_{{\mathbb D}_{R''}} P_{R'}$ for all $R' > R''$, and letting first
$R'$ and then $R''$ go to $\infty$, we see from the discussion above
(with $R \to \infty$ as well) that any two loops are connected by a
finite ``path" of intermediate loops, $P$-a.s. \fbox{} \\

\noindent {\bf Proof of Theorem~\ref{thm-conf-inv}.}
Combining Russo-Seymour-Welsh type arguments for
percolation (see, for example, Lemma~3 of~\cite{ksz})
with Theorem~\ref{main-result1}, we know that $P$-a.s.
there exists a bounded continuum nonsimple loop that
surrounds both $\lambda_1$ and $\lambda_2$, so that
$\tilde D_1$ and $\tilde D_2$ are both bounded.
We can then take $R<\infty$ such that $\lambda_1$ and $\lambda_2$
(and therefore $\tilde D_1$ and $\tilde D_2$) are both contained
in the disc ${\mathbb D}_R$ with probability tending to
$1$ as $R \to \infty$.

Consider now the continuum construction inside the disc
${\mathbb D}_R$ for some large $R$.
Let $\lambda'_1$ (resp., $\lambda'_2$) be the smallest loop
surrounding $D_1$ (resp., $D_2$) produced by the construction
and let $\tilde D'_1$ (resp., $\tilde D'_2$) be the connected
component of ${\mathbb R}^2 \setminus \lambda'_1$ (resp.,
${\mathbb R}^2 \setminus \lambda'_2$) that contains $D_1$
(resp., $D_2$).
It follows from the previous observation and from (the proof of)
Theorem~\ref{thm-therm-lim} that as $R \to \infty$, $\tilde D'_1$
(resp., $\tilde D'_2$) is (with probability tending to $1$)
distributed like $\tilde D_1$ (resp., $\tilde D_2$) and moreover
the loop configuration inside $\tilde D'_1$ (resp., $\tilde D'_2$)
is distributed by $P_{\tilde D_1}$ (resp., $P_{\tilde D_2}$).

This already proves the first claim of the theorem, since it is
clear from the continuum construction inside ${\mathbb D}_R$
that the loop configurations inside $\tilde D'_1$ and $\tilde D'_2$
are independent.
It also means that in order to complete the proof of theorem, it
suffices to prove the second claim for the case of the continuum
construction inside ${\mathbb D}_R$, for all large $R$.
In order to do that, we consider a modified discrete construction
inside ${\mathbb D}_R$, as explained below.
In view of the above observations, we take $R$ large and condition
on the existence inside ${\mathbb D}_R$ of two disjoint loops,
$\lambda^{\delta}_1$ and $\lambda^{\delta}_2$, surrounding $D_1$
and $D_2$ respectively, and let $\tilde D^{\delta}_1$ (resp.,
$\tilde D^{\delta}_2$) be the domain of
${\mathbb D}_R^{\delta} \setminus \Gamma(\lambda^{\delta}_1)$
(resp., ${\mathbb D}_R^{\delta} \setminus \Gamma(\lambda^{\delta}_2)$)
containing $D_1$ (resp., $D_2$).

The modified discrete construction inside ${\mathbb D}_R$ is
analogous to the ``ordinary" one except inside the domains
$\tilde D^{\delta}_1$ and $\tilde D^{\delta}_2$, where the
the exploration paths are coupled to a continuum construction
inside the unit disc in the following way.
Roughly speaking, the discrete construction inside $\tilde D^{\delta}_1$
is one in which the $(x,y)$ pairs (the starting and ending points of
the exploration paths) at each step are chosen to be closest to the
$(\phi_{\delta}(a),\phi_{\delta}(b))$ points in $\tilde D^{\delta}_1$
mapped from the unit disc ${\mathbb D}$ via $\phi_{\delta}$, where the
pairs $(a,b)$ are those that appear at the corresponding steps of the
continuum construction inside ${\mathbb D}$ and $\phi_{\delta}$ is a
certain
conformal map from ${\mathbb D}$ onto $\tilde D^{\delta}_1$,
as specified below.
The discrete construction inside $\tilde D^{\delta}_2$ is coupled
in the same way to the \emph{same} continuum construction inside
$\mathbb D$ via a certain conformal map $\psi_{\delta}$ from ${\mathbb D}$
onto $\tilde D^{\delta}_2$.

The conformal map $\phi_{\delta}$ will be defined for $\delta$
sufficiently small and is specified in the following way.
We fix a point $z_0$ in $\tilde D'_1$ and denote by $\phi$ the
unique conformal map from $\mathbb D$ onto $\tilde D'_1$ such that
$\phi(0)=z_0$ and $\phi'(0)>0$.
For $\delta$ sufficiently small, so that $z_0$ is contained in
$\tilde D^{\delta}_1$, we let $\phi_{\delta}$ be the unique
conformal map from $\mathbb D$ onto $\tilde D^{\delta}_1$ such
that $\phi_{\delta}(0)=z_0$ and $\phi'_{\delta}(0)>0$.

We also denote by $\psi$ the unique conformal map from $\mathbb D$
onto $\tilde D'_2$ such that $\psi(0)=g(z_0)$ and
$\text{sign}(\psi'(0)) = \text{sign}(g'(z_0))$, where $g$ is any
fixed conformal map from $\tilde D'_1$ to $\tilde D'_2$.
Note that, by the uniqueness part of Riemann's mapping theorem,
we can conclude that $\psi = g \circ \phi$.
For $\delta$ sufficiently small, so that $g(z_0)$ is contained
in $\tilde D^{\delta}_2$, we let $\psi_{\delta}$ be the unique
conformal map from $\mathbb D$ onto $\tilde D^{\delta}_2$ such
that $\psi_{\delta}(0)=g(z_0)$ and
$\text{sign}(\psi'_{\delta}(0)) = \text{sign}(g'(z_0))$.

As $\delta \to 0$, $\tilde D^{\delta}_1 \to \tilde D'_1$ and
$\tilde D^{\delta}_2 \to \tilde D'_2$, and by an application
of Rad\'o's theorem (Theorem~\ref{rado-thm}), (the continuous
extensions of) $\phi_{\delta}$ and $\psi_{\delta}$ converge
uniformly in $\overline{\mathbb D}$ to the (continuous extensions
of) $\phi$ and $\psi$ respectively.

We now describe more precisely the modified construction
inside $\tilde D^{\delta}_1$.
Let $\gamma_1$ be the first $SLE_6$ path in $\mathbb D$ from
$a_1$ to $b_1$; because of the conformal invariance of $SLE_6$,
the image $\phi_{\delta}(\gamma_1)$ of $\gamma_1$ under
$\phi_{\delta}$ is a path distributed as the trace of chordal $SLE_6$
in $\tilde D^{\delta}_1$ from $\phi_{\delta}(a_1)$ to $\phi_{\delta}(b_1)$.
The uniform convergence of $\phi_{\delta}$ to $\phi$ and statement (S)
imply that the exploration path $\gamma^{\delta}_1$ inside
$\tilde D^{\delta}_1$ from $x_1$ to $y_1$, chosen to be closest
to $\phi_{\delta}(a_1)$ and $\phi_{\delta}(b_1)$ respectively,
converges in distribution to $\phi(\gamma_1)$, as $\delta \to 0$,
which means that there exists a coupling so that the paths
$\gamma^{\delta}_1$ and $\phi_{\delta}(\gamma_1)$ stay close for
$\delta$ small.

One can use the same strategy as in the
proof of the first part of Theorem~\ref{thm-convergence},
and obtain a discrete construction whose exploration paths
are coupled to the $SLE_6$ paths $\phi_{\delta}(\gamma_k)$
that are the images of the paths $\gamma_k$ in $\mathbb D$.
Then, for this discrete construction inside $\tilde D^{\delta}_1$,
the scaling limit of the exploration paths will be paths
inside $\tilde D'_1$ distributed as the images of the $SLE_6$
paths in $\mathbb D$ under the conformal map $\phi$.

Analogously, for the discrete construction inside $\tilde D^{\delta}_2$,
the scaling limit of the exploration paths will be paths
inside $\tilde D'_2$ distributed as the images of the $SLE_6$
paths in $\mathbb D$ under the (continuous extension of the)
conformal map $\psi: {\mathbb D} \to \tilde D'_2$ which is the
uniform limit of $\psi_{\delta}$ as $\delta \to 0$.

Therefore, the path inside $\tilde D'_2$ obtained as the scaling
limit of an exploration path at a given step of the construction
inside $\tilde D^{\delta}_2$ is the image under the conformal map
$g=\psi \circ \phi^{-1}: \tilde D'_1 \to \tilde D'_2$ of the path
inside $\tilde D'_1$ obtained as the scaling limit of the exploration
path at corresponding step of the construction inside $\tilde D^{\delta}_1$.


In order to conclude the proof, we have to show that
the discrete constructions inside $\tilde D^\delta_1$ and
$\tilde D^\delta_2$ defined above find all the boundaries
(with diameter greater than $\varepsilon$)
in a number of steps that is bounded in probability as
$\delta \to 0$.
This is as in the second part of Theorem~\ref{thm-convergence},
but since we are now dealing with a modified construction
inside general Jordan domains, we need to show that we can
reach the same conclusion.
In order to do so, we will use the fact that the modified
construction is coupled to an ``ordinary" continuum
construction in the unit disc.
We work out the details only for $\tilde D^{\delta}_1$, since
the proof is the same for $\tilde D^{\delta}_2$.

>From the second part of Theorem~\ref{thm-convergence} it
follows that, for any fixed $\varepsilon > 0$ and
$M < \infty$, the probability that the number of steps of the
continuum construction in $\mathbb D$ that are necessary to
ensure that only domains with diameter less than $\varepsilon / M$
are present is larger than $C$, goes to zero as $C \to \infty$.
Since $\phi_{\delta}$ (can be extended to a function that) is
continuous in the compact set $\overline{\mathbb D}$, $\phi_{\delta}$
is uniformly continuous and so we can choose
$M_{\delta}=M(\phi_{\delta}) < \infty$ such that any subdomain
of $\mathbb D$ of diameter at most $\varepsilon / M_{\delta}$
is mapped by $\phi_{\delta}$ to a subdomain of $\tilde D^{\delta}_1$
of diameter at most $\varepsilon$.

Since $\phi_{\delta} \to \phi$, as $\delta \to 0$, where $\phi$ (can
be extended to a function that) is continuous in the compact set
$\overline{\mathbb D}$ and is therefore uniformly continuous, we can
choose $M_0=M(\phi) < \infty$ such that any subdomain of $\mathbb D$ of
diameter at most $\varepsilon / M_0$ is mapped by $\phi_{\delta}$
to a subdomain of $\tilde D'_1$ of diameter at most $\varepsilon$,
and moreover such that $\limsup_{\delta \to 0} M_{\delta} \leq M_0$.

This, combined with the coupling between $SLE_6$ paths and
exploration paths inside $\tilde D^{\delta}$, assures that
the number of steps necessary for the new discrete construction
inside $\tilde D^{\delta}_1$ to find all the loops of diameter
at least $\varepsilon$ is bounded in probability as $\delta \to 0$.

Therefore, the scaling limit, as $\delta \to 0$, of the
modified discrete constructions for $\tilde D^{\delta}_1$ and
$\tilde D^{\delta}_2$ give the measures $P_{\tilde D'_1}$
and $P_{\tilde D'_2}$, and it follows by construction that
$P_{\tilde D'_2} = g * P_{\tilde D'_1}$ for \emph{any} conformal
map $g$ from $\tilde D'_1$ onto $\tilde D'_2$.
Since this is true for all large $R$, by letting $R \to \infty$
we can conclude that $P_{\tilde D_2} = f * P_{\tilde D_1}$.~\fbox{}

\bigskip
\bigskip

\noindent {\bf Acknowledgements.} We are grateful to
Greg Lawler, Oded Schramm and Wendelin Werner for various
interesting and useful conversations and to Stas Smirnov
for communications about a paper in preparation.
F.~C. thanks Wendelin Werner for an invitation to Universit\'e
Paris-Sud 11, and Vincent Beffara and Luiz Renato Fontes for
many helpful discussions.
We thank Michael Aizenman, Lai-Sang Young
and an anonymous referee for comments about
the presentation of our results.
F.~C. and C.~M.~N. acknowledge respectively the kind hospitality
of the Courant Institute and of the Vrije Universiteit Amsterdam
where part of this and related work was done.

\end{document}